\documentclass[12pt,reqno]{amsart}
\usepackage{amsmath,amssymb,amsfonts,amscd,latexsym,amsthm,mathrsfs,verbatim,comment,extarrows,diagbox,xypic,dsfont}
\usepackage{hyperref}
\newtheorem{theorem}{Theorem}
\newtheorem{definition}{Definition}
\newtheorem{corollary}{Corollary}
\newtheorem{lemma}{Lemma}
\newtheorem{remark}{Remark}

\newcommand{\s}[1]{\ensuremath{\scriptstyle{#1}}}

\newcommand{\mo}[2]{\ensuremath{M_{3^{#1}}^{#2}}}
\newcommand{\moo}[3]{\ensuremath{M_{3^{#1}~:{#3}}^{#2}}}

\begin{document}

\title{Classification of the congruence classes of $\mathbf{A}_n^5(n\geq 6)$ with 2-torsion free homology}


\author{Zhongjian Zhu}
\address{College of Mathematics, Physics  and Electronic Information Engineering,\\ Wenzhou University, Zhejiang Wenzhou, {\rm 325035}, China}
\email{zhuzhongjian@amss.ac.cn}

\author{Jianzhong Pan$^{*}$}
\address{Hua Loo-Keng Key Mathematical Laboratory, Institute of Mathematics,\\ Academy of Mathematics and Systems Science,
Chinese Academy of Sciences; \\University of Chinese Academy of Sciences, Beijing, {\rm 100190}, China}
\email{pjz@amss.ac.cn}

\thanks{$^*$Corresponding author.}


\begin{abstract}
In this paper, we classify the congruence classes of  $\mathbf{F}^5_{n(2)}$-polyhedra, i.e. $(n-1)$-connected, at most $(n+5)$-dimensional polyhedra with 2-torsion free homology. The proof relies on the matrix problem technique which was developed in the classification of representations of algebras and  applied to homotopy theory by Baues and Drozd.
\end{abstract}

\maketitle

\section{Introduction}
\label{intro}

Let $\mathbf{A}_n^k(n\geq k+1)$ be the subcategories of the stable homotopy category consisting of $(n-1)$-connected polyhedra, i.e. finite CW-complexes, with dimension at most $n+k$; $\mathbf{F}_n^k$ the full subcategory of $\mathbf{A}_n^k$  consisting of polyhedra with torsion free homology groups; $\mathbf{F}_{n(2)}^k$, the full subcategories of $\mathbf{A}_n^k$ consisting of polyhedra with 2-torsion free homology groups; $\mathbf{F}_{n(2,3)}^k$, the full subcategories of $\mathbf{A}_n^k$ consisting of complexes with 2 and 3 torsion free homology groups. All the above categories are full additive if we consider the wedge of two polyhedra as the coproduct of two objects. The classification problem of an full additive category is to find a complete list of its indecomposable isomorphic classes.

Classifying homotopy types of $\mathbf{A}_n^k$-polyhedra in stable range, i.e. with $n\geq k+1$, is a classical and fundamental task of topology. Such polyhedra have been described for $k\leq 3$ \cite{RefBH}, \cite{RefChang} and for $k\geq 4$, their classification is a wild problem, i.e. contains classification of representations of all finitely generated
algebras over a field \cite{RefBD}. Baues and Drozd also consider the classification problem of $\mathbf{F}_n^k(n\geq k+1)$. For $k\leq 5$, such indecomposable homotopy types have been described in  \cite{RefTF5cels}, \cite{RefTF6cels} or \cite{RefDrMSP}, \cite{RefDrMTS}. For $k=6$, Drozd got tame type classification of congruence classes of homotopy types, and proved that, for $k>6$, this problem is wild in \cite{RefDrFP}.

This is the third of a series of our papers devoted to the stable homotopy types of $\mathbf{A}_n^k$. In our first paper \cite{PZ}, we classify indecomposable homotopy types in  $\mathbf{F}_{n(2,3)}^k$ for $k\leq 6$. Since Drozd shows that the existence of the polyhedra with 2-torsion homology makes the classification problem of $\mathbf{A}_n^k$ be wild \cite{RefDrMTS}, it is natural to classify homotopy types with 2-torsion free homology in $\mathbf{A}_n^k$ for $k=4,5,6$  while for $k>6$ this classification problem is also wild. In our second paper \cite{PZ2}, we get all the indecomposable homotopy types of $\mathbf{F}_{n(2)}^4$ by matrix problem technique and this classification problem is tame, i.e. not wild. And in this paper, we will classify  all the indecomposable congruence classes of homotopy types in $\mathbf{F}_{n(2)}^5$ and find that this classification problem is also tame.  The classification problem of $\mathbf{F}_{n(2)}^6$ is still open because of the highly complexity of the corresponding matrix problem.

 Our main theorem, theorem \ref{maintheorem}, is given at the end of this article.
 Section 2 contains some basic notations and facts about stable homotopy category. Section 3 introduces the main technique: bimodule categories and matrix problem. In this section theorem \ref{theorem3.2} and corollary \ref{corollary3.3} establish a connection between bimodule categories and stable homotopy categories. Section 4 reduces the classification problem of $\mathbf{F}_{n(2)}^5$ to a matrix problem $(\mathscr{A},\mathcal{G})$. In section 5, we 2-localize and 3-localize the matrix problem $(\mathscr{A},\mathcal{G})$ to get $(\mathscr{A}(2),\mathcal{G}(2))$ and $(\mathscr{A}(3),\mathcal{G}(3))$ respectively. We use the known results of indecomposable homotopy types of $\mathbf{F}_{n}^5$ in \cite{RefTF6cels} to get the indecomposable isomorphic classes $\mathrm{ind} \mathscr{A}(2)$ of  $(\mathscr{A}(2),\mathcal{G}(2))$ and use the representation theory of bunch of chains \cite{RefFinitquadra} to get the indecomposable isomorphic classes $\mathrm{ind}\mathscr{A}(3)$ of  $(\mathscr{A}(3),\mathcal{G}(3))$. In section 6,  indecomposable congruence classes of $\mathbf{F}_{n(2)}^5$ are obtained by combining indecomposable isomorphic classes of $\mathscr{A}(2)$ with indecomposable isomorphic classes of $\mathscr{A}(3)$. We give a most complicated example to show the detail  procedures of combination.

\section{Preliminaries}
\label{sec:2}

In this paper, ``Space'' means a based space. Denote by $[X,Y]$ the set of homotopy classes of continuous maps $X\rightarrow Y$ and by $\textsf{CW}$ the homotopy category of polyhedra. The suspension functor  $\Sigma : X\mapsto X[1] ~(X[n]=\Sigma^{n}X)$ defines a natural map
 $[X,Y]\rightarrow [X[n],Y[n]]$. Set $\{X,Y\}$  $=\lim\limits_{n\rightarrow{\infty}}$ $ [X[n],Y[n]]$.
 If $\alpha\in [X[n],Y[n]],~\beta\in [Y[m],Z[m]]$, the class $\beta[n]\cdot \alpha[m]\in [X[m+n],Z[m+n]]$ after stabilization is, by definition, the product $\beta\alpha$ of the classes of $\alpha$ and $\beta$ in $\{X,Z\}$.
 Thus we obtain the stable homotopy category of polyhedra $\textsf{CWS}$. Extending $\textsf{CWS}$ by adding formal negative shifts $X[-n] (n\in \mathds{N})$ of  polyhedra and setting $\{X[-n],Y[-m]\}:=\{X[m],Y[n]\}$, one gets the category $\textsf{S}$ of~\cite{RefCo}, which is a fully additive category, and we denote it by $\textsf{CWS}$ too. By the Freudenthal Theorem(\cite{RefRMS} Theorem.6.26), it is easy to get that the suspension functor
   $$\Sigma: \mathbf{X}^{k}_{n} \longrightarrow \mathbf{X}^{k}_{n+1}~~ (\mathbf{X}^{k}_{n}=\mathbf{A}_n^k, \mathbf{F}_{n}^k,  \mathbf{F}_{n(2)}^k)$$
 is an equivalence for all $n > k+1$ and a full representation equivalence (full, dense and reflecting isomorphisms) for $n= k+1$. We should note that if an additive functor $F: \mathcal{C}\rightarrow \mathcal{D}$  is a full representation equivalence, denoted by $ \mathcal{C}\xlongrightarrow{F \simeq_{rep}}  \mathcal{D}$, then it induces an 1-1 correspondence of indecomposable isomorphic classes of objects of these two additive categories.

 Let $\mathcal{C}$ be an additive category with zero object $\ast$ and biproducts $A\oplus B$ for any objects $A,B\in \mathcal{C}$, where $X\in \mathcal{C}$ means that $X$ is an object of $\mathcal{C}$. $X\in\mathcal{C} $ is decomposable if there is an isomorphism $X\cong A\oplus B$ where $A$ and $B$ are not isomorphic to $\ast$, otherwise $X$ is indecomposable. For example, $X\in \textsf{CW}$ (resp. $\textsf{CWS}$) is indecomposable if $X$ is homotopy equivalent (resp. stable homotopy equivalent ) to $X_{1} \vee X_{2} $ implies one of $X_{1}$ and $X_{2}$ is contractible. A decomposition of $X\in \mathcal{C}$ is an isomorphism
 $$X\cong A_{1}\oplus \cdots \oplus A_{n}, ~~n<\infty, $$
 where $A_{i}$ is indecomposable for $i\in \{1,2,\cdots , n\}$. The classification problem of category $\mathcal{C}$ is to find a complete list of indecomposable isomorphism types in $\mathcal{C}$ and describe the possible decompositions of objects in $\mathcal{C}$.

\begin{definition}\label{definition2.1}
We say two polyhedra $X$ and $Y$ are congruent if there is a polyhedron $Z$ such that $X\vee Z\simeq Y\vee Z$(actually, Z can always be chosen as a wedge of spheres ). We write $X\equiv Y$.
\end{definition}

This is an equivalence relation on the objects of $\textsf{CW}$ ($\textsf{CWS}$). We say the congruence class $[X]$ indecomposable if $X\equiv X_{1}\vee X_{2}$ implies that one of $X_{1},X_{2}$ is contractible. From Theorem 4.26 of \cite{RefCo}, it is easy to know that if $X\equiv Y$ and $X$ is an object of $\mathbf{F}_{n(2)}^{k}$(resp. $\mathbf{F}_{n}^{k}$), then so is $Y$ and vise versa.

\section{Techniques}
\label{sec:3}

\begin{definition}\label{definition3.1}
Let $\mathcal{A}$ and $\mathcal{B}$ be small additive categories. $\mathcal{U}$ is an $\mathcal{A}$-$\mathcal{B}$-bimodule,
i.e. a biadditive functor  $\mathcal{A}^{op}\times \mathcal{B}\rightarrow \textbf{Ab}$, the category of abelian groups.
We define the bimodule category $El(\mathcal{U})$ as follows:
\begin{itemize}
  \item  the set of objects is the disjoint union $\bigcup_{A\in \mathcal{A} ,B\in \mathcal{B}}\mathcal{U}(A,B)$.
  \item A morphism $\alpha\rightarrow\beta$, where $\alpha \in \mathcal{U}(A,B), \beta \in \mathcal{U}(A',B')$ is a pair of morphisms $f:A\rightarrow A' $,
  $g:B\rightarrow B'$ such that $g\alpha=\beta f\in \mathcal{U}(A,B')$ (We write $g\alpha$ instead of $\mathcal{U}(1,g)\alpha$ and $\beta f$ instead of $\mathcal{U}(f,1)\beta$).
\end{itemize}
\end{definition}
Obviously $El(\mathcal{U})$ is an  (full) additive category if so are  $\mathcal{A}$ and $\mathcal{B}$.
\\

Suppose $\mathcal{A}$ and $\mathcal{B}$ are two full subcategories of $\textsf{CW}$ (or $\textsf{CWS}$), then we denote by $\mathcal{A}\dag \mathcal{B}$  the full subcategory of $\textsf{CW}$ (or  $\textsf{CWS}$) consisting of cofibers of  maps  $f:A\rightarrow B$, where $A\in \mathcal{A}, B\in \mathcal{B}$. We also denote by $\mathcal{A}\dag_{m} \mathcal{B}$  the full subcategory of $\mathcal{A}\dag \mathcal{B}$  consisting of cofibers of $f:A\rightarrow B$ such that $H_{m}(f)=0$ and denote by $\Gamma(A,B)$ the subgroup of $\{A,B\}$ consisting of maps $f:A\rightarrow B$ such that $H_{m}(f)=0$, where $A\in \mathcal{A}$, $B\in \mathcal{B}$. From Theorem 3.2 of \cite{PZ2}, we have the following theorem and corollary:

\begin{theorem}\label{theorem3.2}
 Let $\mathcal{A}$ and $\mathcal{B}$ be two full subcategories of $\textsf{CWS}$, suppose that  $\{B, A[1]\}=0$ for all $A\in\mathcal{A}, B\in\mathcal{B}$. Consider  $H: \mathcal{A}^{op}\times \mathcal{B}\rightarrow \textbf{Ab}$, i.e.  $(A,B)\mapsto \{A,B\}$, as an $\mathcal{A}$- $\mathcal{B}$-bimodule. Denote by $\mathcal{I}$ the ideal of category $\mathcal{A}\dag \mathcal{B}$ consisting of morphisms which factor through both $\mathcal{B}$ and $\mathcal{A}[1]$, and by $\mathcal{J}$ the ideal of the category $El(H)$ consisting of morphisms $(\alpha, \beta):f\rightarrow f'$ such that $\beta$ factors through $f'$ and $\alpha$ factors through $f$. Then
 \begin{itemize}
 \item [(1)] the functor   $ C:  El(H)\rightarrow  \mathcal{A}\dag \mathcal{B}$ ($f\mapsto C_{f}$) induces an equivalence  $El(H)/\mathcal{J}\simeq \mathcal{A}\dag \mathcal{B}/\mathcal{I}$.
 \item [(2)]  moreover $\mathcal{I}^{2}=0$, hence the projection $\mathcal{A}\dag\mathcal{B}\rightarrow \mathcal{A}\dag \mathcal{B}/\mathcal{I}$ is a representation equivalence.
 \item [(3)] In particular, let $n<m\leq n+k$ and denote by $\mathcal{\widetilde{B}}$ the full subcategory of $\mathbf{F}_{n(2)}^{k}(n\geq k+1)$ consisting of all $(n-1)$-connected polyhedra of dimension at most $m$ and by $\mathcal{\widetilde{A}}$ the full subcategory of $\mathbf{F}_{n(2)}^{k}(n\geq k+1)$ consisting of all $(m-1)$-connected polyhedra of dimension at most $n+k-1$., then
      $$El(H)/\mathcal{J}\xlongrightarrow{\overline{C}~\simeq}\mathcal{\widetilde{A}}\dag \mathcal{\widetilde{B}}/\mathcal{I}\xlongleftarrow{P\simeq_{rep}}\mathcal{\widetilde{A}}\dag \mathcal{\widetilde{B}}.$$
      gives a natural one-to-one correspondence between isomorphic classes of objects of $El(H)/\mathcal{J}$ and $\mathcal{\widetilde{A}}\dag\mathcal{\widetilde{B}}$.
      $\mathbf{F}_{n(2)}^{k}$ is the full subcategory of $\mathcal{\widetilde{A}}\dag\mathcal{\widetilde{B}}$ consisting of $2$-torsion free polyhedra.
 \end{itemize}
 \end{theorem}

\begin{corollary}\label{corollary3.3}
Under conditions of Theorem \ref{theorem3.2}, let $H_{0}$ be an $\mathcal{A}$-$\mathcal{B}$-subbimodule of $H$ such that $f_1af_2=0$ whenever $a\in El(H_{0})$, $f_{i}\in \{B_{i},A_{i}\} (i=1,2)$. Denote by $\mathcal{A}\dag_{H_{0}}\mathcal{B}$ the full subcategory of $\mathcal{A}\dag\mathcal{B}$ consisting of cofibers of $a\in El(H_{0})$.
$\mathcal{I}_{H_{0}}=Mor(\mathcal{A}\dag_{H_{0}}\mathcal{B})\cap\mathcal{I}$ and $\mathcal{J}_{H_{0}}=Mor(El(H_{0}))\cap\mathcal{J}$. Then we have

 \begin{itemize}
   \item [(1)] $\mathcal{J}^2_{H_0}=\mathcal{I}_{H_{0}}^2=0$;
   \item [(2)] $C: El(H_{0})\xlongrightarrow{P \simeq_{rep}}El(H_{0})/\mathcal{J}_{H_{0}}\xlongrightarrow{\overline{C}~\simeq}\mathcal{A}\dag_{H_{0}} \mathcal{B}/\mathcal{I}_{H_{0}}\xlongleftarrow{P \simeq_{rep}}\mathcal{A}\dag_{H_{0}} \mathcal{B}.$
 \end{itemize}
\end{corollary}

If $H_{0}=\Gamma:\mathcal{A}^{op}\times \mathcal{B}\rightarrow \textbf{Ab}$, then $\mathcal{A}\dag_{H_{0}}\mathcal{B}=\mathcal{A}\dag_{m} \mathcal{B}$.

\textbf{Matrix problem} ~ Let $\mathscr{A}$ be a set of matrices which is closed under finite  direct sums of matrices  and let $\mathcal{G}$ denote the set of admissible transformations on $\mathscr{A}$.
 We say $A\cong B$ in $\mathscr{A}$ if $A$ can be transformed to $B$ by admissible transformations, and we say $A$ is decomposable
 if $A\cong A_{1}\bigoplus A_{2}$ for nontrivial $A_{1}, A_{2}\in \mathscr{A}$. The block matrices $\left( \begin{array}{c} A_{1} \\ 0 \\ \end{array} \right)$ and
$ \left( \begin{array}{cc}  A_{1} &0 \\  \end{array}  \right)$ are also thought to be decomposable. The matrix problem $(\mathscr{A}, \mathcal{G})$,
 or simply $\mathscr{A}$, means to classify the indecomposable isomorphic classes of matrices of $\mathscr{A}$  (denoted by $\mathrm{ind}\mathscr{A}$) under admissible transformations $\mathcal{G}$.
Matrix problem $(\mathscr{A}, \mathcal{G})$ is said to be equivalent to matrix problem $(\mathscr{A}', \mathcal{G}')$ if there is a bijective map
$\varphi: \mathscr{A} \rightarrow  \mathscr{A}' $ such that $ A\cong A' $ in $ \mathscr{A} $ if and only if $ \varphi(A)\cong \varphi ( A')$ in $\mathscr{A}' $
and $\varphi(A_{1}\bigoplus A_{2}) = \varphi(A_{1})\bigoplus \varphi(A_{2})$. It is clear that if two matrix problems are equivalent, then there is a  one-to-one correspondence between their indecomposable isomorphic classes.

\section{The reduction of the classification problem of $\mathbf{F}_{n(2)}^5$~$(n\geq 6)$}
\label{sec:matrixproblem}
In the following context, the tabulations $\begin{tabular}{|ccc|}  \hline * &*  &*  \\  * & * & *\\ \hline  \end{tabular}$ represent the matrices or block matrices.  For any category $\mathcal{C}$, denote by $\mathrm{ind}\mathcal{C}$ the set of indecomposable isomorphic classes of objects of $\mathcal{C}$.

Let $M_{t}^{k}$ be the Moore space $M(\mathbb{Z}/t,k)$ ,  $t,k \in \mathbb{N}_+=\{1,2,\cdots,\}$; $C^{l}_{\eta}$ and $C^{l}_{\eta^{2}}$ be the cofiber of maps  $S^{l-1}\xlongrightarrow{\eta}S^{l-2}$ and  $S^{l-1}\xlongrightarrow{\eta^2}S^{l-3}$ respectively, where $\eta$ is the suspension of the Hopf map.

Take $m=n+3$ in Theorem 3.2(3) to get two full subcategories $\tilde{\mathcal{A}}$ and  $\tilde{\mathcal{B}}$ of $\mathbf{F}_{n(2)}^5$, i.e.  $\tilde{\mathcal{A}}=\mathbf{F}_{n+3(2)}^1$, $\tilde{\mathcal{B}}=\mathbf{F}_{n(2)}^3$. So from \cite{RefBH}
$$\begin{array}{l}
    \mathrm{ind}\tilde{\mathcal{A}}=\{S^{n+3},S^{n+4},M_{p^{r}}^{n+3}~|~\text{prime}~p\neq 2,  r\in \mathbb{N}_+\}; \\
    \mathrm{ind}\tilde{\mathcal{B}}=\{S^{n},S^{n+1},S^{n+2},S^{n+3},C_{\eta}^{n+2}, C_{\eta^{2}}^{n+3},C_{\eta}^{n+3},M_{p^{r}}^{n} ,M_{p^{r}}^{n+1},M_{p^{r}}^{n+2}~|~\text{prime}~p\neq 2, r\in \mathbb{N}_+\}.
  \end{array}
$$
From Lemma 4.4 of \cite{PZ} and Lemma 5.1 of \cite{PZ2}, we get
$$\begin{array}{l}
    \{M_{p^{r}}^{n+3}, X\}=0 ~\text{for}~X\in \mathrm{ind}\tilde{\mathcal{B}}; \text{prime}~ p\neq 2,3;  r\in \mathbb{N}_+;\\
    \{Y, M_{p^{s}}^{n+t}\}=0 ~\text{for}~Y\in \mathrm{ind}\tilde{\mathcal{A}}; t=0,1,2; \text{prime}~ p\neq 2,3; s\in \mathbb{N}_+; \\
    \{Y, M_{3^{s}}^{n+2}\}=0~\text{for}~Y\in \mathrm{ind}\tilde{\mathcal{A}};  s\in \mathbb{N}_+.
  \end{array}
$$
Hence we get that
$$\mathrm{ind}(\tilde{\mathcal{A}}\dag \tilde{\mathcal{B}})=\mathrm{ind}(\mathcal{A}\dag\mathcal{B})\cup \{M_{p^{r}}^{n+4},M_{p^{r}}^{n},M_{p^{r}}^{n+1},M_{q^{r}}^{n+2}~|~\text{prime}~ p\neq 2,3;~\text{prime}~ q\neq 2;  r\in \mathbb{N}_+\},$$
where $\mathcal{A}$ and $\mathcal{B}$ are full subcategories of $\mathcal{\widetilde{A}}$ and $\mathcal{\widetilde{B}}$ respectively such that
  \begin{itemize}
  \item [] $\mathrm{ind}\mathcal{A}=\{S^{n+3},S^{n+4},M_{3^{r}}^{n+3}~|~ r\in \mathbb{N}_+\};$
   \item [] $ \mathrm{ind}\mathcal{B}=\{S^{n},S^{n+1},S^{n+2},S^{n+3},C_{\eta}^{n+2}, C_{\eta^{2}}^{n+3},C_{\eta}^{n+3},M_{3^{s}}^{n},M_{3^{s}}^{n+1} ~|~ s\in \mathbb{N}_+\}.$
  \end{itemize}
  By Theorem \ref{theorem3.2} (3), we get
  \begin{lemma}\label{lemmaindFn25}
Let $\mathcal{A}$ and $\mathcal{B}$ be defined as above. Then

$\mathrm{ind}\mathbf{F}_{n(2)}^5=\{X\in \mathrm{ind}(\mathcal{A}\dag\mathcal{B})$ $|~X$ is $2$-torsion free
$\}$ $\cup$ $\{~M^{n+4}_{p^r}, M^{n}_{p^r},$ $ M^{n+1}_{p^r}$ $ M^{n+2}_{q^r}~|~$ prime $p\neq 2,3$, prime $q\neq 2$ and $r\in \mathbb{N}_+\}.$
\end{lemma}

Let $\Gamma : \mathcal{A}^{op}\times\mathcal{B}\rightarrow \textbf{Ab}$, $\Gamma(A,B)=\{g\in \{A,B\}~|~H_{n+3}(g)=0\}$, defined in section \ref{sec:3}, be a sub-bimodule of $\mathcal{A}$-$\mathcal{B}$-bimodule $H: \mathcal{A}^{op}\times\mathcal{B}\rightarrow \textbf{Ab}$, $H(A,B)=\{A,B\}$.

Since the general case reduces to case with $H_{n+3}$ 3-torsion free and another case similar to it, which will be discussed briefly at the end of the paper, we will focus on the classification of 2-torsion free polyhedra with $H_{n+3}$ 3-torsion free.

By the similar proof as in Lemma 5.4 and Corollary 5.5 of \cite{PZ2}, we have

\begin{lemma}\label{lemmaindA+B2tor}
Let $\mathcal{A}$ and $\mathcal{B}$ be defined as above. Then
  $$\begin{array}{rl}
   & \{X\in \mathrm{ind}(\mathcal{A}\dag\mathcal{B})~|~X~\text{is}~2\text{-torsion free}, H_{n+3}X ~\text{is}~ 3\text{-torsion free}
   \} \\
   =&\{M_{p^{r}}^{n+3}~|~ \text{prime}~ p\neq 2,3, r\in \mathbb{N}_{+}\}\cup \{~C(f)~\text{is indecomposable}~|~f\in El(\Gamma)~\}, \\
  \end{array}$$
  where  $ \{~C(f)~\text{is indecomposable}~|~f\in El(\Gamma)~\}=\mathrm{ind}(\mathcal{A}\dag_{n+3}\mathcal{B})\cong \mathrm{ind} El(\Gamma).$
 \end{lemma}

Compute $\Gamma(A,B)$ for $A\in \mathrm{ind}\mathcal{A}$, $B\in \mathrm{ind}\mathcal{B}$ and write them in matrix forms as in \cite{RefDrFP}. For example
$$\Gamma(\mo{r}{n+3}, \mo{s}{n})=\begin{scriptsize}\begin{tabular}{c|cc|}
   \multicolumn{1}{c}{}&\multicolumn{2}{c}{$\mo{r}{n+3}:\s{ n+3~~n+4}$  }\\
       \cline{2-3}
      $\mo{s}{n}:\s{n}$ & $\mathbb{Z}/3$ &0\\
      \quad ${\scriptstyle{n+1}}$ & ~~~~~0& $\mathbb{Z}/3$ \\
        \cline{2-3}
     \end{tabular} \end{scriptsize}~;\Gamma(\mo{r}{n+3}, \mo{s}{n+1})=\begin{scriptsize}\begin{tabular}{c|cc|}
   \multicolumn{1}{c}{}&\multicolumn{2}{c}{$\mo{r}{n+3}:\s{ n+3~~n+4}$  }\\
       \cline{2-3}
      $\mo{s}{n+1}:\s{n+1}$ &~~~~ 0 &$\mathbb{Z}/3$\\
      \quad ${\scriptstyle{n+2}}$ & ~~~~~0& 0 \\
        \cline{2-3}
     \end{tabular} \end{scriptsize}~;$$
$$\Gamma(\mo{r}{n+3}, C_{\eta}^{n+2})=\begin{scriptsize}\begin{tabular}{c|cc|}
   \multicolumn{1}{c}{}&\multicolumn{2}{c}{$\mo{r}{n+3}:\s{ n+3~~n+4}$  }\\
       \cline{2-3}
      $ C_{\eta}^{n+2}:\s{n}$ &$\mathbb{Z}/3$ &0\\
      \quad ${\scriptstyle{n+2}}$ & ~~~~~0& 0 \\
        \cline{2-3}
     \end{tabular} \end{scriptsize}~;\Gamma(\mo{r}{n+3}, C_{\eta^{2}}^{n+3})=\begin{scriptsize}\begin{tabular}{c|cc|}
   \multicolumn{1}{c}{}&\multicolumn{2}{c}{$\mo{r}{n+3}:\s{ n+3~~n+4}$  }\\
       \cline{2-3}
      $ C_{\eta^{2}}^{n+3}:\s{n}$ &$\mathbb{Z}/3$ &0\\
      \quad ${\scriptstyle{n+3}}$ & ~~~~~0& 0 \\
        \cline{2-3}
     \end{tabular} \end{scriptsize}~.$$
     We can get following $Table ~\Gamma(\mathcal{A},\mathcal{B})$.

     $$\begin{tabular}{|c|c|c|cc|cc|cc|c}
   \multicolumn{10}{c}{$Table~\Gamma(\mathcal{A},\mathcal{B})$}\\
                    \multicolumn{10}{c}{}\\
 \hline
   \diagbox{$\mathcal{B}$}{$\mathcal{A}$} & $S^{n+3}$ & $S^{n+4}$& \multicolumn{2}{c|}{$\mo{}{n+3}:\s{n+3~~n+4}$} &  \multicolumn{2}{c|}{$\mo{2}{n+3}:\s{n+3~~n+4}$} &  \multicolumn{2}{c|}{$\mo{3}{n+3}:\s{n+3~~n+4}$} & $\cdots$ \\
  \hline
 $S^{n}$ &$\mathbb{Z}/24$&0&~~~~~$\mathbb{Z}/3$&~0&~~~~~$\mathbb{Z}/3$&~0&~~~~~$\mathbb{Z}/3$&~0&$\cdots$\\
 \hline
  $S^{n+1}$ &$\mathbb{Z}/2$&$\mathbb{Z}/24$&~~~0&$\mathbb{Z}/3$&~~~0&$\mathbb{Z}/3$&~~~0&$\mathbb{Z}/3$&$\cdots$\\
 \hline
  $S^{n+2}$ &$\mathbb{Z}/2$&$\mathbb{Z}/2$&~~~0&0&~~~0&0&~~~0&0&$\cdots$\\
 \hline
 $S^{n+3}$ &0&$\mathbb{Z}/2$&~~~0&0&~~~0&0&~~~0&0&$\cdots$\\
 \hline
  $C_{\eta}^{n+2}:\s{n}$ &$\mathbb{Z}/12$ &0&~~~~$\mathbb{Z}/3$&0&~~~~$\mathbb{Z}/3$&0&~~~~$\mathbb{Z}/3$&0&$\cdots$\\
  ~~$\s{n+2}$ &0&0&~~~~0&0&~~~~0&0&~~~~0&0&$\cdots$\\
  \hline
  $C_{\eta^2}^{n+3}:\s{n}$ &$\mathbb{Z}/12$ &0&~~~~$\mathbb{Z}/3$&0&~~~~$\mathbb{Z}/3$&0&~~~~$\mathbb{Z}/3$&0&$\cdots$\\
  ~~$\s{n+3}$ &0&0&~~~~0&0&~~~~0&0&~~~~0&0&$\cdots$\\
 \hline
  $C_{\eta}^{n+3}:\s{n+1}$ &0 &$\mathbb{Z}/12$&~~~0&$\mathbb{Z}/3$&~~~0&$\mathbb{Z}/3$&~~~0&$\mathbb{Z}/3$&$\cdots$\\
  ~~$\s{n+3}$ &0&0&~~~0&0&~~~0&0&~~~0&0&$\cdots$\\
 \hline
  $\mo{}{n}:\s{n}$ &$\mathbb{Z}/3$ &0&~~~~~$\mathbb{Z}/3$&0&~~~~~$\mathbb{Z}/3$&0&~~~~~$\mathbb{Z}/3$&0&$\cdots$\\
  ~~$\s{n+1}$ &0&$\mathbb{Z}/3$ &~~~0&$\mathbb{Z}/3$&~~~0&$\mathbb{Z}/3$&~~~0&$\mathbb{Z}/3$&$\cdots$\\
 \hline
  $\mo{2}{n}:\s{n}$ &$\mathbb{Z}/3$ &0&~~~~~$\mathbb{Z}/3$&0&~~~~~$\mathbb{Z}/3$&0&~~~~~$\mathbb{Z}/3$&0&$\cdots$\\
  ~~$\s{n+1}$ &0&$\mathbb{Z}/3$ &~~~0&$\mathbb{Z}/3$&~~~0&$\mathbb{Z}/3$&~~~0&$\mathbb{Z}/3$&$\cdots$\\
 \hline
 $\mo{3}{n}:\s{n}$ &$\mathbb{Z}/3$ &0&~~~~~$\mathbb{Z}/3$&0&~~~~~$\mathbb{Z}/3$&0&~~~~~$\mathbb{Z}/3$&0&$\cdots$\\
  ~~$\s{n+1}$ &0&$\mathbb{Z}/3$ &~~~0&$\mathbb{Z}/3$&~~~0&$\mathbb{Z}/3$&~~~0&$\mathbb{Z}/3$&$\cdots$\\
 \hline
 $\cdots$& $\cdots$& $\cdots$& $\cdots$& $\cdots$& $\cdots$& $\cdots$& $\cdots$& $\cdots$& $\cdots$\\
  \hline
  $\mo{}{n+1}:\s{n+1}$ &0 &$\mathbb{Z}/3$&~~~0&$\mathbb{Z}/3$&~~~0&$\mathbb{Z}/3$&~~~0&$\mathbb{Z}/3$&$\cdots$\\
  ~~~~$\s{n+2}$ &0&0&~~~0&0&~~~0&0&~~~0&0&$\cdots$\\
 \hline
  $\mo{2}{n+1}:\s{n+1}$ &0 &$\mathbb{Z}/3$&~~~0&$\mathbb{Z}/3$&~~~0&$\mathbb{Z}/3$&~~~0&$\mathbb{Z}/3$&$\cdots$\\
  ~~~~$\s{n+2}$ &0&0&~~~0&0&~~~0&0&~~~0&0&$\cdots$\\
 \hline
 $\mo{3}{n+1}:\s{n+1}$ &0 &$\mathbb{Z}/3$&~~~0&$\mathbb{Z}/3$&~~~0&$\mathbb{Z}/3$&~~~0&$\mathbb{Z}/3$&$\cdots$\\
  ~~~~$\s{n+2}$ &0&0&~~~0&0&~~~0&0&~~~0&0&$\cdots$\\
 \hline
 $\cdots$& $\cdots$& $\cdots$& $\cdots$& $\cdots$& $\cdots$& $\cdots$& $\cdots$& $\cdots$& $\cdots$\\
 \multicolumn{10}{c}{}\\
 \end{tabular}
$$
   Objects of $El(\Gamma)$ are represented by block matrices $(\gamma_{ij})$ with finite dimension, where block $\gamma_{ij}$ has entries from the $(ij)$-th cell of $Table~ \Gamma(\mathcal{A},\mathcal{B})$. All these block matrices form a matrix set $\mathscr{A}$. Isomorphisms of $El(\Gamma)$ can be represented by following admissible transformations $\mathcal{G}$ which are obtained by computing $\{A,A'\}$ for $A, A'\in \mathrm{ind} \mathcal{A}$ and $\{B,B'\}$ for $B, B'\in \mathrm{ind} \mathcal{B}$ as in Table 2 and Table 3 on page 9 of \cite{RefDrFP}.
\\

\textbf{ Admissible Transformations} $\mathcal{G}$:

 Let $W_{x}$ (resp. $W^{y}$) denote the $x$-horizontal (resp. $y$-vertical) stripe, where
 $x\in \{S^{n},S^{n+1}$, $S^{n+2}$, $S^{n+3}$, $C_{\eta~~~:n}^{n+2}$, $C_{\eta^{2}~:n}^{n+3}$, $C_{\eta~~~:n+1}^{n+3}$$\}$$\coprod$$\{\moo{s}{n}{n}, \moo{s}{n}{n+1}, \moo{s}{n+1}{n+1}~|~$ $s\in\mathbb{N}_+\}$ (resp. $y\in \{S^{k}~|~k=n+3,n+4\}$ $\coprod$ $\{\moo{r}{n+3}{n+3},$ $ \moo{r}{n+3}{n+4}$ $~|~r\in\mathbb{N}_+\}$). Denote by $W^{y}_{x}$ the block corresponding to $x$-horizontal stripe and $y$-vertical stripe. Let dim $W_{x}=$ the number of rows in $W_{x}$, dim $W^{y}=$ the number of columns in $W^{y}$.

\begin{itemize}
  \item [(a)] Elementary transformations of rows (columns) in each horizontal (vertical) stripe;
   \item[(b)]  $W^{S^{n+3}}<W^{S^{n+4}}$; ~~$W^{\moo{r}{n+3}{n+3}}<W^{\moo{r'}{n+3}{n+3}}<W^{S^{n+3}}$; ~~ $W^{S^{n+4}}<W^{\moo{r'}{n+3}{n+4}}<W^{\moo{r}{n+3}{n+4}}$, where $r<r'$;
   \item[(c)] $W_{S^{n+3}}<W_{S^{n+2}}<W_{S^{n+1}}<W_{S^{n}}$;~~$W_{S^{n}}<W_{C_{\eta~~~:n}^{n+2}}$;~~
    $2W_{C_{\eta~~~:n}^{n+2}}<W_{S^{n}}$;~~

     $6W_{S^{n+2}}<W_{C_{\eta~~~:n}^{n+2}} $;~~$W_{C_{\eta^2~:n}^{n+3}}<W_{C_{\eta~~~:n}^{n+2}}$;~~ $2W_{C_{\eta~~~:n}^{n+2}}<W_{C_{\eta^2~:n}^{n+3}}$;~~$6W_{S^{n+3}}<W_{C_{\eta~~~:n+1}^{n+3}}$;~~

    $2W_{C_{\eta~~~:n+1}^{n+3}}<W_{S^{n+1}}$; ~~$W_{S^{n}}<W_{C_{\eta^2~:n}^{n+3}}$;~~$2W_{C_{\eta^2~:n}^{n+3}}<W_{S^{n}}$;~~

    $W_{S^{n}},W_{C_{\eta^2~:n}^{n+3}},W_{C_{\eta~~~:n}^{n+2}}<W_{\moo{s'}{n}{n}}<W_{\moo{s}{n}{n}}$;  ~~~

    $W_{\moo{s}{n}{n+1}}<W_{\moo{s'}{n}{n+1}}<W_{S^{n+1}}<W_{C_{\eta~~~:n+1}^{n+3}}<W_{\moo{t'}{n+1}{n+1}}<W_{\moo{t}{n+1}{n+1}}$ where $s<s'; t<t'$;
   \item [(d)] $W_{\moo{s}{n}{n}}\thicksim W_{\moo{s}{n}{n+1}}$; $W^{\moo{r}{n+3}{n+3}}\thicksim W^{\moo{r}{n+3}{n+4}}$,
\end{itemize}
\begin{itemize}
  \item [$\diamond$] $x<y$ means scalar multiples of rows (columns) of $x$-stripe can be added to rows (columns) of the $y$-stripe;
   $ax<y $ ($a\in\mathbb{N}_+$) means adding $ak$ times of a row (column) of  $x$-stripe to a row (column) of $y$-stripe is admissible where $k$ is an any nonzero integer;
  \item [$\diamond$] $x\thicksim y$ implies that $dim x=dim y$ and it means that the transformations of the $x$-stripe must be the same as those of the $y$-stripe.
\end{itemize}

  \begin{remark}
    When admissible transformations above are performed on block matrix $\gamma=(\gamma_{ij})$, where block $\gamma_{ij}$ has entries from (ij)-cell of $Table~\Gamma(\mathcal{A},\mathcal{B})$,  we should note that
  \begin{itemize}
    \item [(1)]  If (ij)-cell of $Table~\Gamma(\mathcal{A},\mathcal{B})$ is zero, then $\gamma_{ij}$ keeps being zero after admissible transformations;
    \item [(2)]  Adding $1\in \mathbb{Z}/2$ to an element $a\in \mathbb{Z}/24$ gives $a+12$ in $\mathbb{Z}/24$, since $\eta^{3}$ is $12$ in $\mathbb{Z}/24 = \{S^{n+3},S^{n}\}\cong \{S^{n+4},S^{n+1}\} $.
    \item [(3)] The reason for $6W_{S^{n+3}}<W_{C_{\eta~~~:n+1}^{n+3}}$ and $6W_{S^{n+2}}<W_{C_{\eta~~~:n}^{n+2}}$ above is the Proposition~6 (iii) of \cite{RefUnsold}.
  \end{itemize}
  \end{remark}

So we get a matrix problem $(\mathscr{A}, \mathcal{G})$ above such that $\mathrm{ind} \mathscr{A}\cong \mathrm{ind} El(\Gamma)\cong \mathrm{ind} (\mathcal{A}\dag_{n+3}\mathcal{B})$.

\section{The localization of categories and matrix problems}
\label{sec:localization}
For any additive category $\mathcal{C}$ and prime $p$, we denote by $\mathcal{C}(p)$ the category  with the same set of objects,  but with the sets of morphisms $\mathbb{Z}(p)\otimes \mathcal{C}(X,Y)$, where $\mathbb{Z}(p)=\{m/n~|~m,n\in \mathbb{Z}, p\nmid n\}$. Then we have the natural functor $L(p):\mathcal{C}\rightarrow \mathcal{C}(p)$. We call $\mathcal{C}(p)$ the $p$-localization of the category $\mathcal{C}$.

For $X\in \textsf{CWS}$, we denote by $\mathbf{P}(X)$ the set of all primes dividing the order of one of the stable homotopy groups
$\pi_{k}^{s}(X)$ with $k\leq dim X$. We also use the notation $L(p)X$ to denote the $p$-localization of space $X$.  The following theorem is from the Theorem 1.5 of \cite{RefDrFP}
\begin{theorem}\label{theoremlocalcongu}
For any two polyhedra $X, Y\in \textsf{CWS}$ the following properties are equivalent:
\begin{itemize}
  \item [(1)] $X\equiv Y$;
  \item [(2)] $L(p)X\simeq L(p)Y$ for all prime $p$;
  \item [(3)] $L(p)X\simeq L(p)Y$ for all prime $p\in \mathbf{P}(X\vee Y)$.
\end{itemize}
\end{theorem}

Since $\mathbf{P}(X)\subset \{2,3\}$ for all $X\in \mathcal{A}\dag_{n+3}\mathcal{B}$,  from above theorem, in order to get all indecomposable congruence classes of $\mathcal{A}\dag_{n+3}\mathcal{B}$, we compute $\mathrm{ind}(\mathcal{A}\dag_{n+3}\mathcal{B})(2)$ (i.e. $\mathrm{ind}El(\Gamma)(2)$) and $\mathrm{ind}(\mathcal{A}\dag_{n+3}\mathcal{B})(3)$ (i.e. $\mathrm{ind}El(\Gamma)(3)$), and then try to combine them into indecomposable congruence classes of $\mathcal{A}\dag_{n+3}\mathcal{B}$.

Replace $\mathbb{Z}/24$ by $\mathbb{Z}/8$, $\mathbb{Z}/12$ by $\mathbb{Z}/4$  and $\mathbb{Z}/3$ by $0$ in $Table~\Gamma(\mathcal{A},\mathcal{B})$ and then delete the zero stripes, we get the following table which provide the matrix set $\mathscr{A}(2)$
     $$\begin{tabular}{|c|c|c|}
 \hline
   \diagbox{$\mathcal{B}$}{$\mathcal{A}$} & $S^{n+3}$ & $S^{n+4}$ \\
  \hline
 $S^{n}$ &$\mathbb{Z}/8$&0\\
 \hline
  $S^{n+1}$ &$\mathbb{Z}/2$&$\mathbb{Z}/8$\\
 \hline
  $S^{n+2}$ &$\mathbb{Z}/2$&$\mathbb{Z}/2$\\
 \hline
 $S^{n+3}$ &0&$\mathbb{Z}/2$\\
 \hline
  $C_{\eta}^{n+2}:\s{n}$ &$\mathbb{Z}/4$ &0\\
  ~~$\s{n+2}$ &0&0\\
  \hline
  $C_{\eta^2}^{n+3}:\s{n}$ &$\mathbb{Z}/4$ &0\\
  ~~$\s{n+2}$ &0&0\\
 \hline
  $C_{\eta}^{n+3}:\s{n+1}$ &0 &$\mathbb{Z}/4$\\
  ~~$\s{n+3}$ &0&0\\
 \hline
 \multicolumn{3}{c}{}\\
 \end{tabular}$$

Admissible transformations $\mathcal{G}(2)$  are obtained from both the restriction of $\mathcal{G}$ on the matrix set  $\mathscr{A}(2)$ and the multiplications of a row(column) by an odd number $a$.
Then we get the matrix problem $(\mathscr{A}(2), \mathcal{G}(2))$  such that $\mathrm{ind}\mathscr{A}(2)\cong \mathrm{ind}(\mathcal{A}\dag_{n+3}\mathcal{B})(2)$. Note that the matrix problem is just the same as the ``2-localization" of the matrix problem $(\mathscr{A}^{5}, \mathcal{G}^{5})$ in \cite{PZ}, which satisfies $\mathrm{ind}\mathscr{A}^{5}\cong \mathrm{ind}\mathbf{F}_{n}^{5}$. By finding the corresponding matrix forms of indecomposable $\mathbf{F}_{n}^{5}$-polyhedra \cite{RefTF6cels} in the matrix set $\mathscr{A}^{5}$ and then 2-localizing them, we get the following list of $\mathrm{ind}\mathscr{A}(2)$.

 \textbf{List} $\mathrm{ind}\mathscr{A}(2)$ :

   $S^n; S^{n+1}; S^{n+2}; S^{n+3}; S^{n+4}; S^{n+5}$;~~$\s{C_{\upsilon}^{n+4}=}\begin{tabular}{c|c|}
   \multicolumn{1}{c}{}&\multicolumn{1}{c} {$\s{S^{n+3}}$}\\
       \cline{2-2}
     $\s{S^{n}}$ & $\s \upsilon$ \\
       \cline{2-2}
     \end{tabular}$;~~ $\s{C_{\omega}^{n+5}=}\begin{tabular}{c|c|}
   \multicolumn{1}{c}{}&\multicolumn{1}{c} {$\s{S^{n+4}}$}\\
       \cline{2-2}
     $\s{S^{n+1}}$ & $\s \omega$ \\
       \cline{2-2}
     \end{tabular}$;~~

$\s{C_{\eta}^{n+2}=}\begin{tabular}{c|c|}
   \multicolumn{1}{c}{}&\multicolumn{1}{c} {$\s{S^{n+1}}$}\\
       \cline{2-2}
     $\s{S^{n}}$ & $\s \eta$ \\
       \cline{2-2}
     \end{tabular}$;~~ $\s{C_{\eta}^{n+3}=}\begin{tabular}{c|c|}
   \multicolumn{1}{c}{}&\multicolumn{1}{c} {$\s{S^{n+2}}$}\\
       \cline{2-2}
     $\s{S^{n+1}}$ & $\s \eta$ \\
       \cline{2-2}
     \end{tabular}$;~~ $\s{C_{\eta}^{n+4}=}\begin{tabular}{c|c|}
   \multicolumn{1}{c}{}&\multicolumn{1}{c} {$\s{S^{n+3}}$}\\
       \cline{2-2}
     $\s{S^{n+2}}$ & $\s \eta$ \\
       \cline{2-2}
     \end{tabular}$;~~ $\s{C_{\eta}^{n+5}=}\begin{tabular}{c|c|}
   \multicolumn{1}{c}{}&\multicolumn{1}{c} {$\s{S^{n+4}}$}\\
       \cline{2-2}
     $\s{S^{n+3}}$ & $\s \eta$ \\
       \cline{2-2}
     \end{tabular}$;~~

     $\s{C_{\eta^2}^{n+3}=}\begin{tabular}{c|c|}
   \multicolumn{1}{c}{}&\multicolumn{1}{c} {$\s{S^{n+2}}$}\\
       \cline{2-2}
     $\s{S^{n}}$ & $\s \eta^2$ \\
       \cline{2-2}
     \end{tabular}$;~~
   $\s{C_{\eta^2}^{n+4}=}\begin{tabular}{c|c|}
   \multicolumn{1}{c}{}&\multicolumn{1}{c} {$\s{S^{n+3}}$}\\
       \cline{2-2}
     $\s{S^{n+1}}$ & $\s \eta^2$ \\
       \cline{2-2}
     \end{tabular}$;~~
     $\s{C_{\eta^2}^{n+5}=}\begin{tabular}{c|c|}
   \multicolumn{1}{c}{}&\multicolumn{1}{c} {$\s{S^{n+4}}$}\\
       \cline{2-2}
     $\s{S^{n+2}}$ & $\s \eta^2$ \\
       \cline{2-2}
     \end{tabular}$;~~

    $\s{(\eta\upsilon\eta)_{0}^1=\begin{tabular}{c|c|}
   \multicolumn{1}{c}{}&\multicolumn{1}{c} {$\s{S^{n+3}}$}\\
       \cline{2-2}
     $\s{S^{n+2}}$ & $\s \eta$ \\
      \cline{2-2}
          $\s{C_{\eta}^{n+2}}:\s{n}$ & $\s{\upsilon}$   \\
      \quad ${\s{n+2}}$ & \s 0 \\
        \cline{2-2}
     \end{tabular}}$~;~~  $\s{(\eta\omega\eta)_{1}^1=\begin{tabular}{c|c|}
   \multicolumn{1}{c}{}&\multicolumn{1}{c} {$\s{S^{n+4}}$}\\
       \cline{2-2}
     $\s{S^{n+3}}$ & $\s \eta$ \\
      \cline{2-2}
          $\s{C_{\eta}^{n+3}}:\s{n+1}$ & $\s{\omega}$   \\
      \quad ${\s{n+3}}$ & \s 0 \\
        \cline{2-2}
     \end{tabular}}$~;~~ $\s{(\eta^2\upsilon\eta^2)_{0}^1=\begin{tabular}{c|c|}
   \multicolumn{1}{c}{}&\multicolumn{1}{c} {$\s{S^{n+3}}$}\\
       \cline{2-2}
     $\s{S^{n+1}}$ & $\s \eta^2$ \\
      \cline{2-2}
          $\s{C_{\eta^2}^{n+3}}:\s{n}$ & $\s{\upsilon}$   \\
      \quad ${\s{n+3}}$ & \s 0 \\
        \cline{2-2}
     \end{tabular}}$~;~~

    $\s{(\eta^2\upsilon\eta)_{0}^1=\begin{tabular}{c|c|}
   \multicolumn{1}{c}{}&\multicolumn{1}{c} {$\s{S^{n+3}}$}\\
       \cline{2-2}
     $\s{S^{n+2}}$ & $\s \eta$ \\
      \cline{2-2}
          $\s{C_{\eta^2}^{n+3}}:\s{n}$ & $\s{\upsilon}$   \\
      \quad ${\s{n+3}}$ & \s 0 \\
        \cline{2-2}
     \end{tabular}}$~;~~   $\s{(\eta\upsilon\eta^2)_{0}^1=\begin{tabular}{c|c|}
   \multicolumn{1}{c}{}&\multicolumn{1}{c} {$\s{S^{n+3}}$}\\
       \cline{2-2}
     $\s{S^{n+1}}$ & $\s \eta^2$ \\
      \cline{2-2}
          $\s{C_{\eta}^{n+2}}:\s{n}$ & $\s{\upsilon}$   \\
      \quad ${\s{n+2}}$ & \s 0 \\
        \cline{2-2}
     \end{tabular}}$~;~~ $\s{(\eta\omega\eta^2)_{1}^1=\begin{tabular}{c|c|}
   \multicolumn{1}{c}{}&\multicolumn{1}{c} {$\s{S^{n+4}}$}\\
       \cline{2-2}
     $\s{S^{n+2}}$ & $\s \eta^2$ \\
      \cline{2-2}
          $\s{C_{\eta}^{n+3}}:\s{n+1}$ & $\s{\omega}$   \\
      \quad ${\s{n+3}}$ & \s 0 \\
        \cline{2-2}
     \end{tabular}}$~;~~

    $\s{(\upsilon\eta^2)_{0}^0=\begin{tabular}{c|c|}
   \multicolumn{1}{c}{}&\multicolumn{1}{c} {$\s{S^{n+3}}$}\\
       \cline{2-2}
     $\s{S^{n}}$ & $\s \upsilon$ \\
      \cline{2-2}
          $\s{S^{n+1}}$ & $\s{\eta^2}$   \\
      \cline{2-2}
     \end{tabular}}$~;~~ $\s{(\omega\eta^2)_{1}^0=\begin{tabular}{c|c|}
   \multicolumn{1}{c}{}&\multicolumn{1}{c} {$\s{S^{n+4}}$}\\
       \cline{2-2}
     $\s{S^{n+1}}$ & $\s \omega$ \\
      \cline{2-2}
          $\s{S^{n+2}}$ & $\s{\eta^2}$   \\
      \cline{2-2}
     \end{tabular}}$~;~~  $\s{(\upsilon\eta)_{0}^0=\begin{tabular}{c|c|}
   \multicolumn{1}{c}{}&\multicolumn{1}{c} {$\s{S^{n+3}}$}\\
       \cline{2-2}
     $\s{S^{n}}$ & $\s \upsilon$ \\
      \cline{2-2}
          $\s{S^{n+2}}$ & $\s{\eta}$   \\
      \cline{2-2}
     \end{tabular}}$~;~~

    $\s{(\omega\eta)_{1}^0=\begin{tabular}{c|c|}
   \multicolumn{1}{c}{}&\multicolumn{1}{c} {$\s{S^{n+4}}$}\\
       \cline{2-2}
     $\s{S^{n+1}}$ & $\s \omega$ \\
      \cline{2-2}
          $\s{S^{n+3}}$ & $\s{\eta}$   \\
      \cline{2-2}
     \end{tabular}}$~;~~    $\s{(\eta\upsilon)_{0}^1=\begin{tabular}{c|c|}
   \multicolumn{1}{c}{}&\multicolumn{1}{c} {$\s{S^{n+3}}$}\\
       \cline{2-2}
    $\s{C_{\eta}^{n+2}}:\s{n}$ & $\s \upsilon$\\
      \quad ${\s{n+2}}$ & \s 0 \\
      \cline{2-2}
     \end{tabular}}$~;~~ $\s{(\eta^2\omega\eta^2)_{1}^1=\begin{tabular}{c|c|c|}
   \multicolumn{1}{c}{}&\multicolumn{1}{c} {$\s{S^{n+3}}$}&\multicolumn{1}{c} {$\s{S^{n+4}}$}\\
       \cline{2-3}
   $\s{S^{n+1}}$ &$\s \eta^2$  & $\s \omega$\\
      \cline{2-3}
     $\s{S^{n+2}}$ & $\s 0$ &  $\s \eta^2$\\
      \cline{2-3}
     \end{tabular}}$~;

     $\s{(\eta^2\upsilon)_{0}^1=\begin{tabular}{c|c|}
   \multicolumn{1}{c}{}&\multicolumn{1}{c} {$\s{S^{n+3}}$}\\
       \cline{2-2}
    $\s{C_{\eta^2}^{n+3}}:\s{n}$ & $\s \upsilon$\\
      \quad ${\s{n+3}}$ & \s 0 \\
      \cline{2-2}
     \end{tabular}}$~;~~ $\s{(\eta\omega)_{1}^1=\begin{tabular}{c|c|}
   \multicolumn{1}{c}{}&\multicolumn{1}{c} {$\s{S^{n+4}}$}\\
       \cline{2-2}
    $\s{C_{\eta}^{n+3}}:\s{n+1}$ & $\s \omega$\\
      \quad ${\s{n+3}}$ & \s 0 \\
      \cline{2-2}
     \end{tabular}}$~;~~ $\s{(\upsilon\eta^2\omega)_{0}^0=\begin{tabular}{c|c|c|}
   \multicolumn{1}{c}{}&\multicolumn{1}{c} {$\s{S^{n+3}}$}&\multicolumn{1}{c} {$\s{S^{n+4}}$}\\
       \cline{2-3}
   $\s{S^{n}}$ & $\s \upsilon$ & \s 0\\
      \cline{2-3}
     $\s{S^{n+1}}$ & $\s \eta^2$ & $\s \omega$ \\
      \cline{2-3}
     \end{tabular}}$~;~~

        $\s{(\eta^2\omega\eta)_{1}^1=\begin{tabular}{c|c|c|}
   \multicolumn{1}{c}{}&\multicolumn{1}{c} {$\s{S^{n+3}}$}&\multicolumn{1}{c} {$\s{S^{n+4}}$}\\
       \cline{2-3}
   $\s{S^{n+1}}$ &$\s \eta^2$  & $\s \omega$\\
      \cline{2-3}
     $\s{S^{n+3}}$ & $\s 0$ &  $\s \eta$\\
      \cline{2-3}
     \end{tabular}}$~;~~$\s{(\upsilon\eta^2\omega\eta)_{0}^0=\begin{tabular}{c|c|c|}
   \multicolumn{1}{c}{}&\multicolumn{1}{c} {$\s{S^{n+3}}$}&\multicolumn{1}{c} {$\s{S^{n+4}}$}\\
       \cline{2-3}
   $\s{S^{n}}$ &$\s \upsilon$  & $\s 0$\\
      \cline{2-3}
   $\s{S^{n+1}}$ & $\s \eta^2$ &  $\s \omega$\\
      \cline{2-3}
    $\s{S^{n+3}}$ & $\s 0$ &  $\s \eta$\\
      \cline{2-3}
     \end{tabular}}$~;~~

     $\s{(\upsilon\eta^2\omega\eta^2)_{0}^0=\begin{tabular}{c|c|c|}
   \multicolumn{1}{c}{}&\multicolumn{1}{c} {$\s{S^{n+3}}$}&\multicolumn{1}{c} {$\s{S^{n+4}}$}\\
       \cline{2-3}
   $\s{S^{n}}$ &$\s \upsilon$  & $\s 0$\\
      \cline{2-3}
   $\s{S^{n+1}}$ & $\s \eta^2$ &  $\s \omega$\\
      \cline{2-3}
    $\s{S^{n+2}}$ & $\s 0$ &  $\s \eta^2$\\
      \cline{2-3}
     \end{tabular}}$~;~~$\s{(\eta\upsilon\eta^2\omega)_{0}^1=\begin{tabular}{c|c|c|}
   \multicolumn{1}{c}{}&\multicolumn{1}{c} {$\s{S^{n+3}}$}&\multicolumn{1}{c} {$\s{S^{n+4}}$}\\
       \cline{2-3}
   $\s{C_{\eta}^{n+2}}:\s{n}$ & $\s \upsilon$ & \s 0\\
      \quad ${\s{n+2}}$ & \s 0 & \s 0\\
      \cline{2-3}
    $\s{S^{n+1}}$ &$\s \eta^2$ & $\s \omega$ \\
      \cline{2-3}
     \end{tabular}}$~;~~

  $\s{(\eta^2\upsilon\eta^2\omega)_{0}^1=\begin{tabular}{c|c|c|}
   \multicolumn{1}{c}{}&\multicolumn{1}{c} {$\s{S^{n+3}}$}&\multicolumn{1}{c} {$\s{S^{n+4}}$}\\
       \cline{2-3}
   $\s{C_{\eta^2}^{n+3}}:\s{n}$ & $\s \upsilon$ & \s 0\\
      \quad ${\s{n+3}}$ & \s 0 & \s 0\\
      \cline{2-3}
    $\s{S^{n+1}}$ &$\s \eta^2$ & $\s \omega$ \\
      \cline{2-3}
     \end{tabular}}$~;~~$\s{(\eta^2\omega)_1^1}=\begin{tabular}{c|c|c|}
   \multicolumn{1}{c}{}&\multicolumn{1}{c} {$\s{S^{n+3}}$}&\multicolumn{1}{c} {$\s{S^{n+4}}$}\\
       \cline{2-3}
     $\s{S^{n+1}}$ & $\s \eta^2$& $\s \omega$ \\
       \cline{2-3}
     \end{tabular}$;~~

  $\s{(\eta\upsilon\eta^2\omega\eta)_{0}^1=\begin{tabular}{c|c|c|}
   \multicolumn{1}{c}{}&\multicolumn{1}{c} {$\s{S^{n+3}}$}&\multicolumn{1}{c} {$\s{S^{n+4}}$}\\
       \cline{2-3}
   $\s{C_{\eta}^{n+2}}:\s{n}$ & $\s \upsilon$ & \s 0\\
      \quad ${\s{n+2}}$ & \s 0 & \s 0\\
      \cline{2-3}
    $\s{S^{n+1}}$ &$\s \eta^2$ & $\s \omega$ \\
      \cline{2-3}
    $\s{S^{n+3}}$ &$\s 0$ & $\s \eta$ \\
      \cline{2-3}
     \end{tabular}}$~;~~ $\s{(\eta^2\upsilon\eta^2\omega\eta)_{0}^1=\begin{tabular}{c|c|c|}
   \multicolumn{1}{c}{}&\multicolumn{1}{c} {$\s{S^{n+3}}$}&\multicolumn{1}{c} {$\s{S^{n+4}}$}\\
       \cline{2-3}
   $\s{C_{\eta^2}^{n+3}}:\s{n}$ & $\s \upsilon$ & \s 0\\
      \quad ${\s{n+3}}$ & \s 0 & \s 0\\
      \cline{2-3}
    $\s{S^{n+1}}$ &$\s \eta^2$ & $\s \omega$ \\
      \cline{2-3}
    $\s{S^{n+3}}$ &$\s 0$ & $\s \eta$ \\
      \cline{2-3}
     \end{tabular}}$~;~~

 $\s{(\eta\upsilon\eta^2\omega\eta^2)_{0}^1=\begin{tabular}{c|c|c|}
   \multicolumn{1}{c}{}&\multicolumn{1}{c} {$\s{S^{n+3}}$}&\multicolumn{1}{c} {$\s{S^{n+4}}$}\\
       \cline{2-3}
   $\s{C_{\eta}^{n+2}}:\s{n}$ & $\s \upsilon$ & \s 0\\
      \quad ${\s{n+2}}$ & \s 0 & \s 0\\
      \cline{2-3}
    $\s{S^{n+1}}$ &$\s \eta^2$ & $\s \omega$ \\
      \cline{2-3}
    $\s{S^{n+2}}$ &$\s 0$ & $\s \eta^2$ \\
      \cline{2-3}
     \end{tabular}}$~;~~$\s{(\eta^2\upsilon\eta^2\omega\eta^2)_{0}^1=\begin{tabular}{c|c|c|}
   \multicolumn{1}{c}{}&\multicolumn{1}{c} {$\s{S^{n+3}}$}&\multicolumn{1}{c} {$\s{S^{n+4}}$}\\
       \cline{2-3}
   $\s{C_{\eta^2}^{n+3}}:\s{n}$ & $\s \upsilon$ & \s 0\\
      \quad ${\s{n+3}}$ & \s 0 & \s 0\\
      \cline{2-3}
    $\s{S^{n+1}}$ &$\s \eta^2$ & $\s \omega$ \\
      \cline{2-3}
    $\s{S^{n+2}}$ &$\s 0$ & $\s \eta^2$ \\
      \cline{2-3}
     \end{tabular}}$
 \\
 \\
 where $\upsilon~(\omega) \in\{1,2,4\}\subset\mathbb{Z}/8$ above with row label $S^n$ ($S^{n+1}$);   $\upsilon~(\omega) \in\{1,2\}\subset\mathbb{Z}/4$ above with row label $C_{\eta~~~:n}^{n+2}, C_{\eta^{2}~:n}^{n+3}$ ($C_{\eta~:n+1}^{n+3}$).

Similarly, replace $\mathbb{Z}/24$, $\mathbb{Z}/12$ by $\mathbb{Z}/3$, $\mathbb{Z}/2$  by 0 in $Table~\Gamma(\mathcal{A},\mathcal{B})$, then we get the  matrix problem  $(\mathscr{A}(3),  \mathcal{G}(3))$ which satisfies
 $\mathrm{ind} \mathscr{A}(3)\cong \mathrm{ind} (\mathcal{A}\dag_{n+3}\mathcal{B})(3)$.
Note that in the matrix problem $(\mathscr{A}(3), \mathcal{G}(3))$, any row-transformation among  $S^n$, $C_{\eta~~~:n}^{n+2}$, $C_{\eta^{2}~:n}^{n+3}$ horizonal stripes is admissible, so we can combine these three horizonal stripes into one stripe, labeled by $e_{0}$-stripe. Similarly,  we also can combine $S^{n+1}$-stripe and $C_{\eta~:n+1}^{n+3}$-stripe into one stripe, labeled by $e'_{0}$-stripe. We should remark here that in the remainder of the paper, the label $e_{0}$ (resp. $e'_{0}$) sometimes also means an element of the set $\{S^n, C_{\eta~~~:n}^{n+2}, C_{\eta^{2}~:n}^{n+3}\}$ (resp. the set $\{S^{n+1}, C_{\eta~:n+1}^{n+3}\}$).
In order to simplify writing we also denote $\moo{s}{n}{n}$-stripe by $e_{s}$;  $\moo{s}{n}{n+1}$-stripe by $\tilde{e}_{s}$;  $\moo{s}{n+1}{n+1}$-stripe by $e'_{s}$; $\moo{r}{n+3}{n+3}$ by $f_{r}$; $\moo{r}{n+3}{n+4}$ by $\tilde{f}_{r}$; $S^{n+3}$-stripe by $f_{0}$ and $S^{n+4}$-stripe by $f'_{0}$.

 Now the matrix problem $(\mathscr{A}(3), \mathcal{G}(3))$ becomes the following example of a bunch of chains (cf. \cite{RefBunchChains} or Appendix of \cite{RefFinitquadra}. We use the notations of the latter paper ):
\begin{eqnarray}
&&\mathfrak{E}_1=\{e_{0}<e_{s'}<e_{s},~~  s<s'\},~~~~\mathfrak{F}_1=\{f_r<f_{r'}<f_{0},~~r<r'\}     \nonumber\\
&&\mathfrak{E}_2=\{\tilde{e}_{s}<\tilde{e}_{s'}<e'_{0}<e'_{t'}<e'_{t},~~  s<s',t<t'\},~~~~\mathfrak{F}_2=\{f'_{0}<\tilde{f}_{r'}<\tilde{f}_{r},~~r<r'\}     \nonumber
\end{eqnarray}
with the following equivalence relation $\thicksim$
$$e_{s}\thicksim \tilde{e}_{s}, s\geq 1,~~~~ f_{r}\thicksim \tilde{f}_{r}, r\geq 1. $$

Let $\mathfrak{E}=\mathfrak{E}_{1}\coprod \mathfrak{E}_{2}$; $\mathfrak{F}=\mathfrak{F}_{1}\coprod \mathfrak{F}_{2}$; $\mathfrak{X}_i=\mathfrak{E}_i\coprod \mathfrak{F}_i$ ($i=1,2$); $\mathfrak{X}= \coprod_{i}\mathfrak{X}_i$.
For $x\in\mathfrak{X}$,  denote $[x]$ the cardinal number of the set $\{y\in\mathfrak{X} |y\neq x,  x\thicksim y \}$, hence $[x]\in\{0,1\}$. We also write $x-y$ if $x\in \mathfrak{E}_i$, $y\in\mathfrak{F}_i $ or vice versa (for some $i\in\{1,2\}$).
We define a word which is a sequence $x_1r_1x_2r_2\cdots x_{l-1}r_{l-1}x_{l}$ with $x_{i}\in \mathfrak{X}$ and $r_{i}\in \{\thicksim, -\}$ such that (i) $r_{i}\neq r_{i+1}$; (2) $x_{i}r_{i}x_{i+1}(1\leq i<l)$ according to the definition of the relations $\thicksim$ and $-$ given above; (3) if $r_1=-(r_l=-)$, then $[x_1]=0$ ($[x_l]=0$).

All the indecomposable isomorphic classes for the above bunch of chains are described as ``string objects" and ``band objects''.

A string object is represented by a word defined above, and two string objects are isomorphic if and only if  the two words are inverse to each other (the inverse of the word $x_1r_1x_2r_2\cdots x_{l-1}r_{l-1}x_{l}$  is $x_lr_{l-1}x_{l-1}r_{l-2}\cdots r_2x_2r_{1}x_{l}$ ). In the following we list all the string objects and their corresponding matrix forms in the matrix problem $(\mathscr{A}(3),  \mathcal{G}(3))$:

\textbf{List} $\mathrm{ind}\mathscr{A}(3)$ :

  \begin{itemize}
     \item Type 1 : $x_1,x_l\in \mathfrak{E}$

    (i)~$[x_1]=0, r_1=-, [x_l]=0, r_{l-1}=- $, then $l=2t$,

     $x_1-x_2\thicksim x_3-x_4\thicksim\cdots x_{2t-2}\thicksim  x_{2t-1}-x_{2t}$

    \begin{tabular}{c|cc|ccc|cc|}
 \multicolumn{1}{c}{}&\multicolumn{1}{c}{$x_{2}$}&\multicolumn{1}{c}{$x_{3}$}&\multicolumn{1}{c}{$x_{6}$}&\multicolumn{1}{c}{$\cdots$}&
 \multicolumn{1}{c}{$x_{l-3}$}&\multicolumn{1}{c}{$x_{l-2}$}&\multicolumn{1}{c}{$x_{l-1}$}\\
 \cline{2-8}
    $x_{1}$& 1& 0 & 0& ... &0&0& ~~0\\
    \cline{2-8}
    $x_{4}$&0 & 1 & 0& ... &0& 0& ~~0\\
    $x_{5}$& 0 & 0 & 1 & ...&0& 0 &~~0 \\
      \cline{2-8}
    ...& ...& ... &... &... &... &...& ...\\
      \cline{2-8}
      $x_{l-4}$& 0 & 0 &0 &...& 1& 0& ~~0 \\
    $x_{l-3}$& 0 & 0 &0 &...& 0& 1& ~~0 \\
      \cline{2-8}
   $x_{l}$& 0 &0 & 0 &... & 0 &0 &~~1 \\
    \cline{2-8}
 \end{tabular}

\bigskip
     (ii)~$[x_1]=1, r_{1}=\thicksim, [x_l]=0, r_{l-1}=-$, then $l=2t+1$,

     $x_1\thicksim x_2-x_3\thicksim x_4-\cdots x_{2t-1}\thicksim  x_{2t}-x_{2t+1}$

       \begin{tabular}{c|cc|ccc|cc|}
 \multicolumn{1}{c}{}&\multicolumn{1}{c}{$x_{3}$}&\multicolumn{1}{c}{$x_{4}$}&\multicolumn{1}{c}{$x_{7}$}&\multicolumn{1}{c}{$\cdots$}&
 \multicolumn{1}{c}{$x_{l-3}$}&\multicolumn{1}{c}{$x_{l-2}$}&\multicolumn{1}{c}{$x_{l-1}$}\\
 \cline{2-8}
    $x_{1}$& 0& 0 & 0& ... &0&0& ~~0\\
    $x_{2}$&1 & 0 & 0& ... &0& 0& ~~0\\
     \cline{2-8}
    $x_{5}$& 0 & 1 & 0 & ...&0& 0 &~~0 \\
    ...& ...& ... &... &... &... &...& ...\\
      \cline{2-8}
      $x_{l-4}$& 0 & 0 &0 &...& 1& 0& ~~0 \\
    $x_{l-3}$& 0 & 0 &0 &...& 0& 1& ~~0 \\
      \cline{2-8}
   $x_{l}$& 0 &0 & 0 &... & 0 &0 &~~1 \\
    \cline{2-8}
 \end{tabular}

    \bigskip

 (iii)~$[x_1]=1, r_{1}=\thicksim, [x_l]=1, r_{l-1}=\thicksim$, then $l=2t$,

  $x_1\thicksim x_2-x_3\thicksim x_4-\cdots x_{2t-2}- x_{2t-1}\thicksim x_{2t}$

       \begin{tabular}{c|cc|cc|c|cc|}
 \multicolumn{1}{c}{}&\multicolumn{1}{c}{$x_{3}$}&\multicolumn{1}{c}{$x_{4}$}&\multicolumn{1}{c}{$x_{7}$}&\multicolumn{1}{c}{$x_{8}$}&
 \multicolumn{1}{c}{$\cdots$}&\multicolumn{1}{c}{$x_{l-4}$}&\multicolumn{1}{c}{$x_{l-3}$}\\
 \cline{2-8}
    $x_{1}$& 0& 0 & 0& 0 &...&0& ~~0\\
    $x_{2}$&1 & 0 & 0& 0 &...& 0& ~~0\\
     \cline{2-8}
    $x_{5}$& 0 & 1 & 0 &0&...& 0 &~~0 \\
      $x_{6}$& 0 & 0 &1 &0& ...& 0& ~~0 \\
          \cline{2-8}
         ...& ...& ... &... &0 &... &...& ...\\
        \cline{2-8}
    $x_{l-1}$& 0 & 0 &0 &0& ...& 0& ~~1 \\
   $x_{l}$& 0 &0 & 0 &0 & ... &0 &~~0 \\
    \cline{2-8}
 \end{tabular}

    \bigskip

     \item Type 2 :  $x_1,x_l\in \mathfrak{F}$

    (i)~$[x_1]=0, r_1=-, [x_l]=0, r_{l-1}=- $, then $l=2t$,

     $x_1-x_2\thicksim x_3-x_4\thicksim x_5 \cdots x_{2t-2}\thicksim  x_{2t-1}-x_{2t}$

    \begin{tabular}{c|c|cc|c|cc|c|}
 \multicolumn{1}{c}{}&\multicolumn{1}{c}{$x_{1}$}&\multicolumn{1}{c}{$x_{4}$}&\multicolumn{1}{c}{$x_{5}$}&\multicolumn{1}{c}{$\cdots$}&
 \multicolumn{1}{c}{$x_{l-4}$}&\multicolumn{1}{c}{$x_{l-3}$}&\multicolumn{1}{c}{$x_{l}$}\\
 \cline{2-8}
    $x_{2}$& 1& 0 & 0& ... &0&0& ~~0\\
    $x_{3}$&0 & 1 & 0& ... &0& 0& ~~0\\
     \cline{2-8}
    $x_{5}$& 0 & 0 & 1 & ...&0& 0 &~~0 \\
    ...& ...& ... &... &... &... &...& ...\\
      \cline{2-8}
    $x_{l-2}$& 0 & 0 &0 &...& 0& 1& ~~0 \\
   $x_{l-1}$& 0 &0 & 0 &... & 0 &0 &~~1 \\
    \cline{2-8}
 \end{tabular}

    \bigskip

     (ii)~$[x_1]=1, r_1=\thicksim, [x_l]=0, r_{l-1}=- $, then $l=2t+1$,

     $x_1\thicksim x_2-x_3\thicksim x_4 \cdots x_{2t-1}\thicksim  x_{2t}-x_{2t+1}$

    \begin{tabular}{c|cc|cc|cc|c|}
 \multicolumn{1}{c}{}&\multicolumn{1}{c}{$x_{1}$}&\multicolumn{1}{c}{$x_{2}$}&\multicolumn{1}{c}{$x_{5}$}&\multicolumn{1}{c}{$\cdots$}&
 \multicolumn{1}{c}{$x_{l-4}$}&\multicolumn{1}{c}{$x_{l-3}$}&\multicolumn{1}{c}{$x_{l}$}\\
 \cline{2-8}
    $x_{3}$& 0& 1 & 0& ... &0&0& ~~0\\
    $x_{4}$&0 & 0 & 1& ... &0& 0& ~~0\\
     \cline{2-8}
    $x_{7}$& 0 & 0 & 0 & ...&0& 0 &~~0 \\
    ...& ...& ... &... &... &... &...& ...\\
      \cline{2-8}
    $x_{l-2}$& 0 & 0 &0 &...& 0& 1& ~~0 \\
   $x_{l-1}$& 0 &0 & 0 &... & 0 &0 &~~1 \\
    \cline{2-8}
 \end{tabular}

\bigskip

 (iii)~$[x_1]=1, r_{1}=\thicksim, [x_l]=1, r_{l-1}=\thicksim$, then $l=2t$,

  $x_1\thicksim x_2-x_3\thicksim x_4-\cdots x_{2t-2}- x_{2t-1}\thicksim x_{2t}$

       \begin{tabular}{c|cc|cc|c|cc|}
 \multicolumn{1}{c}{}&\multicolumn{1}{c}{$x_{1}$}&\multicolumn{1}{c}{$x_{2}$}&\multicolumn{1}{c}{$x_{5}$}&\multicolumn{1}{c}{$x_{6}$}&
 \multicolumn{1}{c}{$\cdots$}&\multicolumn{1}{c}{$x_{l-1}$}&\multicolumn{1}{c}{$x_{l}$}\\
 \cline{2-8}
    $x_{3}$& 0& 1 & 0& 0 &...&0& ~~0\\
    $x_{4}$&0 & 0 & 1& 0 &...& 0& ~~0\\
     \cline{2-8}
    $x_{7}$& 0 & 0 & 0 &1&...& 0 &~~0 \\
      $x_{8}$& 0 & 0 &0 &0& ...& 0& ~~0 \\
          \cline{2-8}
         ...& ...& ... &... &0 &... &...& ...\\
        \cline{2-8}
    $x_{l-4}$& 0 & 0 &0 &0& ...& 1& ~~0 \\
   $x_{l-3}$& 0 &0 & 0 &0 & ... &0 &~~0 \\
    \cline{2-8}
 \end{tabular}

    \bigskip

     \item Type 3 : $x_1\in \mathfrak{E},  x_l\in \mathfrak{F}$ ( $x_1\in \mathfrak{F},  x_l\in \mathfrak{E}$)

    (i)~$[x_1]=0, r_1=-, [x_l]=0, r_{l-1}=- $, then $l=2t$,

     $x_1-x_2\thicksim x_3-x_4\thicksim x_5 \cdots x_{2t-2}\thicksim  x_{2t-1}-x_{2t}$

    \begin{tabular}{c|cc|cc|cc|c|}
 \multicolumn{1}{c}{}&\multicolumn{1}{c}{$x_{2}$}&\multicolumn{1}{c}{$x_{3}$}&\multicolumn{1}{c}{$x_{6}$}&\multicolumn{1}{c}{$\cdots$}&
 \multicolumn{1}{c}{$x_{l-4}$}&\multicolumn{1}{c}{$x_{l-3}$}&\multicolumn{1}{c}{$x_{l}$}\\
 \cline{2-8}
    $x_{1}$& 1& 0 & 0& ... &0&0& ~~0\\
    \cline{2-8}
    $x_{4}$&0 & 1 & 0& ... &0& 0& ~~0\\
    $x_{5}$& 0 & 0 & 1 & ...&0& 0 &~~0 \\
     \cline{2-8}
    ...& ...& ... &... &... &... &...& ...\\
      \cline{2-8}
    $x_{l-2}$& 0 & 0 &0 &...& 0& 1& ~~0 \\
   $x_{l-1}$& 0 &0 & 0 &... & 0 &0 &~~1 \\
    \cline{2-8}
 \end{tabular}

 \bigskip

      (ii)~$[x_1]=1, r_1=\thicksim, [x_l]=0, r_{l-1}=- $, then $l=2t+1$,

     $x_1\thicksim x_2-x_3\thicksim x_4 \cdots x_{2t-1}\thicksim  x_{2t}-x_{2t+1}$

    \begin{tabular}{c|cc|cc|cc|c|}
 \multicolumn{1}{c}{}&\multicolumn{1}{c}{$x_{3}$}&\multicolumn{1}{c}{$x_{4}$}&\multicolumn{1}{c}{$x_{7}$}&\multicolumn{1}{c}{$\cdots$}&
 \multicolumn{1}{c}{$x_{l-4}$}&\multicolumn{1}{c}{$x_{l-3}$}&\multicolumn{1}{c}{$x_{l}$}\\
 \cline{2-8}
    $x_{1}$& 0& 0 & 0& ... &0&0& ~~0\\
    $x_{2}$&1 & 0 & 0& ... &0& 0& ~~0\\
     \cline{2-8}
    $x_{5}$& 0 & 1 & 0 & ...&0& 0 &~~0 \\
     $x_{6}$& 0 & 0& 1 & ...&0& 0 &~~0 \\
     \cline{2-8}
    ...& ...& ... &... &... &... &...& ...\\
      \cline{2-8}
    $x_{l-2}$& 0 & 0 &0 &...& 0& 1& ~~0 \\
   $x_{l-1}$& 0 &0 & 0 &... & 0 &0 &~~1 \\
    \cline{2-8}
 \end{tabular}

   \bigskip

           (iii)~$[x_1]=0, r_1=-, [x_l]=1, r_{l-1}=\thicksim$, then $l=2t+1$,

     $x_1-x_2\thicksim x_3- x_4\thicksim x_5- \cdots x_{2t-2}\thicksim x_{2t-1}- x_{2t}\thicksim x_{2t+1}$

    \begin{tabular}{c|cc|cc|cc|}
 \multicolumn{1}{c}{}&\multicolumn{1}{c}{$x_{2}$}&\multicolumn{1}{c}{$x_{3}$}&\multicolumn{1}{c}{$x_{6}$}&
 \multicolumn{1}{c}{$\cdots$}&\multicolumn{1}{c}{$x_{l-1}$}&\multicolumn{1}{c}{$x_{l}$}\\
 \cline{2-7}
    $x_{1}$& 1& 0 & 0& ... &0& ~~0\\
      \cline{2-7}
    $x_{4}$&0 & 1 & 0& ... &0&  ~~0\\
    $x_{5}$& 0 &0 &1 & ...&0&~~0 \\
     \cline{2-7}
     $x_{8}$& 0 & 0& 0 & ...&0&~~0 \\
    ...& ...& ... &... &... &... & ...\\
      \cline{2-7}
    $x_{l-3}$& 0 & 0 &0 &...& 0&  ~~0 \\
   $x_{l-2}$& 0 &0 & 0 &... & 1 &~~0 \\
    \cline{2-7}
 \end{tabular}

   \bigskip

        (iv)~$[x_1]=1, r_1=\thicksim, [x_l]=1, r_{l-1}=\thicksim$, then $l=2t$,

     $x_1\thicksim x_2-x_3\thicksim x_4-  \cdots x_{2t-3}\thicksim x_{2t-2}- x_{2t-1}\thicksim x_{2t}$

    \begin{tabular}{c|cc|cc|c|cc|}
 \multicolumn{1}{c}{}&\multicolumn{1}{c}{$x_{3}$}&\multicolumn{1}{c}{$x_{4}$}&\multicolumn{1}{c}{$x_{7}$}&\multicolumn{1}{c}{$x_{8}$}&
 \multicolumn{1}{c}{$\cdots$}&\multicolumn{1}{c}{$x_{l-1}$}&\multicolumn{1}{c}{$x_{l}$}\\
 \cline{2-8}
    $x_{1}$&0& 0 & 0&  0&... &0& ~~0\\
    $x_{2}$&1 & 0 & 0& 0& ... &0&  ~~0\\
    \cline{2-8}
    $x_{5}$& 0 &1 &0 &  0&...&0&~~0 \\
     $x_{6}$& 0 & 0& 1 & 0& ...&0&~~0 \\
       \cline{2-8}
    ...& ...& ... &... &... &... &... & ...\\
      \cline{2-8}
    $x_{l-4}$& 0 & 0 &0 & 0&...& 0&  ~~0 \\
   $x_{l-3}$& 0 &0 & 0 & 0&... & 1 &~~0 \\
    \cline{2-8}
 \end{tabular}

   \end{itemize}

 We give some examples of string objects.
 \begin{itemize}
   \item [] Type 1(i)  $e_0-f_1\thicksim \tilde{f}_1-\tilde{e}_3\thicksim e_3-f_4\thicksim \tilde{f}_4-e'_t ~ (t=0~ or~ t\geq 1 )$

    \begin{tabular}{c|cc|cc|}
 \multicolumn{1}{c}{}&\multicolumn{1}{c}{$\s{\mo{1}{n+3}}:\s{n+3}$}&\multicolumn{1}{c}{$\s{n+4}$}&\multicolumn{1}{c}{$\s{\mo{4}{n+3}}:\s{n+3}$}&\multicolumn{1}{c}{$\s n+4$}\\
 \cline{2-5}
       $\s{S^n}$,$\s{C_{\eta~~~:n}^{n+2}}$,$\s{C_{\eta^{2}~:n}^{n+3}}$& 1& 0 &0& ~~0\\
    \cline{2-5}
    $\s{\mo{3}{n}}:\s{n+1}$&0 & 1 &0& ~~0\\
    $\s{n}$& 0 & 0 & 1&~~0 \\
      \cline{2-5}
   $\s{S^{n+1},C_{\eta}^{n+3},~\text{or}~\mo{t}{n+1}}:\s{n+1}$& 0 &0 & 0  &~~1 \\
    \cline{2-5}
 \end{tabular}

   \item []Type 2(ii) $\tilde{f}_3\thicksim f_3-e_3\thicksim \tilde{e}_3- \tilde{f}_4\thicksim f_4- e_5\thicksim \tilde{e}_5- \tilde{f}_3\thicksim f_3$

     \begin{tabular}{c|cc|cc|cc|}
 \multicolumn{1}{c}{}&\multicolumn{1}{c}{$\s{\mo{3}{n+3}:\s{n+4}}$}&\multicolumn{1}{c}{$\s{n+3}$}&\multicolumn{1}{c}{$\s{\mo{4}{n+3}}:\s{n+4}$}&\multicolumn{1}{c}{$\s n+3$}&\multicolumn{1}{c}{$\s{\mo{3}{n+3}:\s{n+4}}$}&\multicolumn{1}{c}{$\s{n+3}$}\\
    \cline{2-7}
    $\s{\mo{3}{n}}:\s{n}$&0 & 1 &0& 0 &0& ~~0\\
    $\s{n+1}$& 0 & 0 & 1 & 0 &0&~~0 \\
      \cline{2-7}
       $\s{\mo{5}{n}}:\s{n}$&0 & 0 &0& 1 &0& ~~0\\
    $\s{n+1}$& 0 & 0 & 0 & 0 &1&~~0 \\
      \cline{2-7}
 \end{tabular}

   \item []Type 3(i)   $e'_t-\tilde{f}_3\thicksim f_3- e_3\thicksim \tilde{e}_3-\tilde{f}_4\thicksim f_4- e_5\thicksim \tilde{e}_5-f'_0      ~(t=0~or ~t\geq 1) $

        \begin{tabular}{c|cc|cc|c|}
 \multicolumn{1}{c}{}&\multicolumn{1}{c}{$\s{\mo{4}{n+3}:\s{n+4}}$}& \multicolumn{1}{c}{$\s n+3$}&\multicolumn{1}{c}{$\s{\mo{4}{n+3}:\s{n+4}}$}&\multicolumn{1}{c}{$\s n+3$}& \multicolumn{1}{c}{$\s{S^{n+4}}$}\\
    \cline{2-6}
    $\s{S^{n+1},C_{\eta}^{n+3},~\text{or}~\mo{t}{n+1}}:\s{n+1}$&1& 0 &0&0& ~~0\\
    \cline{2-6}
    $\s{\mo{3}{n}:\s{n}}$&0& 1 &0& 0&~~0\\
    $\s{n+1}$& 0 &0 & 1 &0&~~0 \\
      \cline{2-6}
       $\s{\mo{5}{n}:\s{n}}$&0 & 0 &0&1& ~~0\\
    $\s{n+1}$& 0 & 0 & 0 &0&~~1\\
      \cline{2-6}
 \end{tabular}

 \end{itemize}

  \bigskip

 We call a word $w=x_1\thicksim x_2-x_3\thicksim x_4-\cdots x_{4m-3}\thicksim x_{4m-2}- x_{4m-1}\thicksim x_{4m} ~\text{with}~ x_1-x_{4m}$  non-periodic cycle if it satisfies $w\neq w^{[k]}$ for $0<k<l$, where $w^{[k]}=x_{k+1}r_{k+1}\cdots r_{k-1}x_{k}$. A band object, denoted by $B(w,z,\pi)$, is represented by a non-periodic cycle $w$, a positive integer $z$ and a unital irreducible polynomial $\pi\neq t$ of degree $v$ from $\mathbb{Z}/3[t]$. The following is its matrix form in $(\mathscr{A}(3),  \mathcal{G}(3))$:

    \begin{tabular}{c|cc|cc|c|cc|}
 \multicolumn{1}{c}{}&\multicolumn{1}{c}{$x_{3}$}&\multicolumn{1}{c}{$x_{4}$}&\multicolumn{1}{c}{$x_{7}$}&\multicolumn{1}{c}{$x_{8}$}&
 \multicolumn{1}{c}{$\cdots$}&\multicolumn{1}{c}{$x_{l-1}$}&\multicolumn{1}{c}{$x_{l}$}\\
 \cline{2-8}
    $x_{1}$&0& 0 & 0&  0&... &0& ~~$F$ \\
    $x_{2}$&$I$ & 0 & 0& 0& ... &0&  ~~0\\
    \cline{2-8}
    $x_{5}$& 0 &$I$  &0 &  0&...&0&~~0 \\
     $x_{6}$& 0 & 0&$I$  & 0& ...&0&~~0 \\
       \cline{2-8}
    ...& ...& ... &... &... &... &... & ...\\
      \cline{2-8}
    $x_{l-4}$& 0 & 0 &0 & 0&...& 0&  ~~0 \\
   $x_{l-3}$& 0 &0 & 0 & 0&... & $I$  &~~0 \\
    \cline{2-8}
 \end{tabular}~~~~~~~~~~~~~
\newline
\newline
where $I$ is the identity matrix and $F$ is the Frobenius matrix with the characteristic polynomial $\pi$. All blocks are of size $zv\times zv$.  The isomorphisms between above band objects are described in Appendix of \cite{RefFinitquadra}. Hence we get the complete list of $\mathrm{ind}\mathscr{A}(3)$ above.
     \bigskip

\section{Combining $\mathrm{ind}\mathscr{A}(2)$ and $\mathrm{ind}\mathscr{A}(3)$ into congruence classes}
\label{Combine}
\textbf{Some notations}

In this section, the notation $\{A,B\}$ means a set with two elements $A,B$ in it.
\begin{itemize}
  \item  Note that  maps $ \mathbb{Z}/8\times \mathbb{Z}/3 \xlongrightarrow{T_{24}}  \mathbb{Z}/24~,(a ,b )\mapsto 9a+16b$ and  $ \mathbb{Z}/4\times \mathbb{Z}/3 \xlongrightarrow{T_{12}}  \mathbb{Z}/12~,(a,b)\mapsto 9a+4b $ are the inverse of the ring isomorphisms $\mathbb{Z}/24  \xlongrightarrow{L_{24}}  \mathbb{Z}/8\times \mathbb{Z}/3 ~,1\mapsto (1,1)$ and
$\mathbb{Z}/12 \xlongrightarrow{T_{12}} \mathbb{Z}/4\times \mathbb{Z}/3 ~,1\mapsto (1,1)$ respectively. We denote both $T_{24}(u,v)$ and $T_{12}(u,v)$ by $u[v]$.

  \item Call a stripe labeled by Moore space $Moore$-stripe, such as  $\moo{s}{n}{n+1}$-horizonal-stripe, $\moo{r}{n+3}{n+4}$-vertical-stripe, $\cdots$.

  \item  Let string object $M\in \mathrm{ind}\mathscr{A}(3)$.  Let $\mathbb{L}= e_0\coprod e'_0\coprod \{S^{k}~|~k=n+3,n+4\}$.
  \begin{itemize}
    \item Denote $M$ by $M(x)$ when $M$ contains an $x$-row (or column ) $x\in\mathbb{L}$.
    \item Denote $M$ by $M(x,y)$ or $M(y,x)$ when $M$ contains both $x$-row or $x$-column  and $y$-row or $y$-column, $x,y\in \mathbb{L}$.
    \item if $M$ contains only one $x$-row (column) with $x\in \mathbb{L}$ and contains no $y$-column (row) with $y\in \mathbb{L}$,  then we call $M$ $e$-type ($f$-type), denoted by $M(x)_e$ ($M(x)_f$); If $M$ contains nether $x$-row with $x\in e_0\coprod e'_0 $ nor $y$-column with $y\in \{S^{k}~|~k=n+3,n+4\}$, then we call $M$ $m$-type.
    \item  $M(e_0)$ ($M(e'_0)$) means some element in $\{M(x)|x\in e_0\}$ ($\{M(x)|x\in e'_0\}$);$M(e_0)_e$ ($M(e'_0)_e$) means some element in $\{M(x)_e|x\in e_0\}$ ($\{M(x)_e|x\in e'_0\}$);$M(e_0, e'_0)$ or $M(e'_0, e_0)$ means some element in $\{M(x,y)~|~x\in e_0, y\in e'_0\}$; $M(e_0, S^{n+3})$ or $M(S^{n+3},e_0)$ means some element in $\{M(x,S^{n+3})~|~x\in e_0\}$;
        $M(e'_0, S^{n+4})$ or $M(S^{n+4},e'_0)$ means some element in $\{M(y,S^{n+4})~|~y\in e'_0\}$;

  \end{itemize}

  The example Type 1(i) above is $M(e_0, e'_0)$ for $t=0$ and  $M(e_0)_e$ for $t\geq 1$; example Type 2(ii) is a $m$-type; example  Type 3(i) is $M(e'_0,S^{n+4})$ for $t=0$ and $M(S^{n+4})_f$ for $t\geq 1$.

  \item If there is a nonzero entry of matrix  $X\in  \mathrm{ind} \mathscr{A}(2)$ and a nonzero entry of  string object (matrix form)  $M(x)\in \mathrm{ind} \mathscr{A}(3)$ in the same $x$-row or $x$-column, then we call $M(x)$  and $X$ are connected by $x$ or $M(x)$ is connected to $X$ by $x$, denoted by $M(x)-X$ or $X-M(x)$;

  \item $M(x,y)\doteqdot X$ or $X\doteqdot M(x,y)$ means $M(x,y)$ and $X$ are connected by both $x$ and $y$. The following object with matrix form can be represented by $$\s{M_1(S^{n+3}, S^{n+4})\doteqdot (\upsilon\eta^2\omega)_0^0<\begin{array}{lc}
                                          \s M_2(S^n)_e &\\
                                         \s M_3(S^{n+1},S^{n+4})-C_{\eta^2}^{n+5}&
                                        \end{array}}$$

$\begin{tabular}{c|c|c|cc|cc|c|}
 \multicolumn{1}{c}{}&\multicolumn{1}{c}{$\s{S^{n+3}}$}&\multicolumn{1}{c}{$\s S^{n+4}$}&\multicolumn{1}{c}{$\s{\mo{4}{n+3}}:\s{n+3}$}&\multicolumn{1}{c}{$\s n+4$}&\multicolumn{1}{c}{$\s{\mo{5}{n+3}}:\s{n+3}$}&\multicolumn{1}{c}{$\s n+4$}&\multicolumn{1}{c}{$\s{S^{n+4}}$}\\
 \cline{2-8}
   $\s S^n$& $\s\upsilon$& \s 0 & \s 1&\s 0 &\s 0& ~~\s 0&\s 0\\
      \cline{2-8}
   $\s S^{n+1}$& $\s\eta^2$& $\s \omega$ &\s  0& \s 0 &\s 0& ~~\s 1&\s 0\\
      \cline{2-8}
   $\s{\mo{3}{n}}:\s{n}$&\s 1 & \s 0& \s 0&\s 0&\s 0&  ~~\s 0&\s 0\\
    $\s{n+1}$& \s 0 &\s 1 &\s 0 & \s 0&\s 0&~~\s 0 &\s 0\\
     \cline{2-8}
      $\s{\mo{5}{n}}:\s{n}$&\s 0 & \s 0 & \s 0&\s 0&\s 1&  ~~\s 0&\s 0\\
    $\s{n+1}$& \s 0 &\s 0 &\s 0 & \s 0&\s 0&~~\s 0&\s 1 \\
    \cline{2-8}
     $\s{S^{n+2}}$&\s  0 &\s 0 &\s 0 & \s 0&\s 0&~~\s 0& $\s\eta^2$\\
    \cline{2-8}
 \end{tabular}$~;

\item If $X_1$, $X_2$, $\cdots$ , $X_l$ ($X_i\in \mathrm{ind}\mathscr{A}(2)\coprod  \mathrm{ind}\mathscr{A}(3)$ for $i=1,2,\cdots,l$) are connected together, then we also call the combination, denoted by $X$, ``string" or ``string object". If we remove some objects $X_{i_1}$, $X_{i_2}$, $\cdots$,  in the string $X$, then we get a sub-string of $X$. For example, for the above string,  $\s{M_1(S^3,S^4)\doteqdot (\upsilon\eta^2\omega)_0^0- \s M_2(S^n)_e}$ is one of its sub-string by removing $M_3(S^{n+1},S^{n+4})$ and $C_{\eta^2}^{n+5}$;
    We call a string  indecomposable if the matrix corresponding to the string is indecomposable.
\end{itemize}

\textbf{Method to combine}

For a polyhedron $X\in \mathcal{A}\dag_{n+3}\mathcal{B}$, let $A_X$ be a matrix in $\mathscr{A}$ realizing $X$ and let
$A_{X}(p)$  be the corresponding matrices of $L(p)X$  in  $\mathscr{A}(p)$ $(p=2,3)$.  We are able to transform $A_{X}(p)$ to $A'_{X}(p)$ which is a direct sum of some matrices listed in $\mathrm{ind}\mathscr{A}(p)$ by admissible transformations. Then we combine
$A'_{X}(2)$  and $A'_{X}(3)$ into a matrix $A'_{X}\in \mathscr{A}$ in the following way:

 let $w$, $u$ and $v$  be the (i,j)-th entry of $A_X$, $A'_{X}(2)$ and $A'_{X}(3)$  respectively, then
 $$ \text{ (i,j)-th entry of} ~A'_{X}=\left\{
                                           \begin{array}{ll}
                                             u[v], & \hbox{if $w\in \mathbb{Z}/24$ or $\mathbb{Z}/12$ ;} \\
                                             u, & \hbox{if $w\in \mathbb{Z}/2$ (implies $v$=0) ;} \\
                                             v, & \hbox{if $w\in \mathbb{Z}/3$ (implies $u$=0) .}
                                           \end{array}
                                         \right.$$
From Theorem \ref{theoremlocalcongu}, we know that the polyhedra represented by matrix $A'_X$ is congruent to $X$.

\textbf{ Key points to split indecomposable congruent classes from a combined matrix}
\begin{itemize}
  \item [(I)] If the band object $B(w,z,\pi)$ or string object $M$ which is m-type is a summand of $A'_{X}(3)$, then $B(w,z,\pi)$ ($M$) is also a summand of $A'_X$, since $B(w,z,\pi)$ ($M$)is a matrix with only $Moore$-stripes and all entries in $Moore$-stripes of $A'_{X}(2)$ are zero. So in the remainder of this section we just consider the case that $A'_{X}(3)$  is a direct sum of string objects which are not $m$-types.
  \item  [(II)] $A'_{X}(3)$ above is a matrix with every row and every column having at most one nonzero entry 1.
  \item [(III)] Note that for a string object $M\in \mathrm{ind}\mathscr{A}(3)$, then $M$ is one of $m$-type, $e$-type and $f$-type or $M\in\{M(e_0,e'_0), M(e_0,S^{n+3}), M(e'_0,S^{n+4}), M(S^{n+3}, S^{n+4})\}$. If $M$ is $e$-type or $f$-type, then it can be connected to at most one polyhedron in  $\mathrm{ind} \mathscr{A}(2)$; if $M\in\{M(e_0,e'_0), M(e_0,S^{n+3})$, $M(e'_0,S^{n+4})$, $M(S^{n+3}, S^{n+4})\}$, then it can be connected to at most two polyhedra in  $\mathrm{ind} \mathscr{A}(2)$.

\end{itemize}

Let $\s{Z(e_0,\upsilon,\omega)\in\left\{
                        \begin{array}{ll}
                         \s{ \{(\eta\upsilon\eta^2\omega\eta)_{0}^1, (\eta\upsilon\eta^2\omega)_{0}^{1}, (\eta\upsilon\eta^2\omega\eta^2)_{0}^1\}}, & \hbox{$\s{e_0=C_{\eta~~~:n}^{n+2}}$;} \\
                          \s{\{ (\upsilon\eta^2\omega\eta)_{0}^0,(\upsilon\eta^2\omega)_{0}^0, (\upsilon\eta^2\omega\eta^2)_{0}^0\}}, & \hbox{$\s{e_0=S^n}$;} \\
                         \s{ \{(\eta^2\upsilon\eta^2\omega\eta)_{0}^1,(\eta^2\upsilon\eta^2\omega)_{0}^1,(\eta^2\upsilon\eta^2\omega\eta^2)_{0}^1 \} }, & \hbox{$\s{e_0=C_{\eta^{2}~:n}^{n+3}}$.}
                        \end{array}
                      \right.};$

$\s{Z(e_0,S^{n+1},\upsilon)=\left\{
                                        \begin{array}{ll}
                                         \s{ (\eta\upsilon\eta^2)_0^1}, & \hbox{$\s{e_0=C_{\eta~~~:n}^{n+2}}$;} \\
                                         \s{ (\upsilon\eta^2)_0^0}, & \hbox{$\s{e_0=S^{n}}$;} \\
                                         \s{  (\eta^2\upsilon\eta^2)_0^1}, & \hbox{$\s{e_0=C_{\eta^{2}~:n}^{n+3}}$.}
                                        \end{array}
                                      \right.}; $~~
$\s{Z(e_0,\upsilon)\in\left\{
                                        \begin{array}{ll}
                                        \s{ \{ (\eta\upsilon)_0^1,(\eta\upsilon\eta)_0^1 \}}, & \hbox{$ \s{e_0=C_{\eta~~~:n}^{n+2}}$;} \\
                                          \s{ \{ C_{\upsilon}^{n+4},(\upsilon\eta)_0^0 \}}, & \hbox{$ \s{e_0=S^{n}}$;} \\
                                           \s{  \{ (\eta^2\upsilon)_0^1,(\eta^2\upsilon\eta)_0^1 \}}, & \hbox{$ \s{e_0=C_{\eta^{2}~:n}^{n+3}}$.}
                                        \end{array}
                                      \right.}; $

                      $\s{Z(e'_0,\omega)\in\left\{
                        \begin{array}{ll}
                         \s{ \{C_{\omega}^{n+5}, (\omega\eta^2)_{1}^{0}, (\omega\eta)_{1}^0\}}, & \hbox{$ \s{e'_0=S^{n+1}}$;} \\
                          \s{ \{(\eta\omega)_{1}^1,(\eta\omega\eta)_{1}^1,(\eta\omega\eta^2)_{1}^1 \}} , & \hbox{$ \s{e'_0=C_{\eta~:n+1}^{n+3}}$.}
                        \end{array}
                      \right.};$~~$\s{Z(\omega)\in\{(\eta^2\omega)^1_1, (\eta^2\omega\eta)^1_1, (\eta^2\omega\eta^2)^1_1 \}}$

In order to simplify writing,  let $S^{k}$ denote $S^{n+k}$, $C_{\eta}^k$ denote $C_{\eta}^{n+k}$, $C_{\eta^2}^k$ denote $C_{\eta^2}^{n+k}$ and let $C^4\in \{C_{\eta^2}^4, C_{\eta}^4\},~~C^5\in \{C_{\eta^2}^5, C_{\eta}^5\}$.

Now we are able to combine string objects in $\mathrm{ind}\mathscr{A}(3)$ with objects in $\mathrm{ind}\mathscr{A}(2)$ to get the following list of  indecomposable congruence classes.

\textbf{List*}
\begin{itemize}
  \item [$\s{(1)}$] $\s Z(e_{0},\upsilon,\omega)$;
  \item [$\s{(2)}$] $\s{ Z(e_{0},\upsilon_1,\omega_0)-M}, ~~~M\in\{M(S^{1})_{e}, M(S^{1},e_0),M(S^{4})_{f},M(S^{4}, S^{3}) \} $;
  \item [$\s{(3)}$]  $\s{Z(e_{0},\upsilon_1,\omega_0)\doteqdot M(S^{1},S^{4})}$;
  \item [$\s{(4)}$]  $\s{M_{1}-Z(e_{0},\upsilon_1,\omega_0)-M_{2},~~~(M_{1},M_{2})\in\{(M_{1}(S^{1})_{e},M_{2}(S^{4})_{f}), ~(M_{1}(S^{1})_{e},\tiny{M_{2}(S^{4},S^{3}))},~(M_{1}(S^{1},e_{0}),M_{2}(S^{4})_f),}$

      $\s{(M_{1}(e_{0})_{e},M_{2}(S^{3})_{f}),~ (M_{1}(e_{0})_{e},M_{2}(S^{3}, S^4)), ~(M_{1}(S^{3},e_{0}),M_{2}(S^{3})_{f}) \}}$;
 \item [$\s{(5)}$]  $\s{M_{1}(S^1)_{e}-Z(e_{0},\upsilon_1,\omega_0)-M_{2}(S^4,S^3)-C^4}$
 and its sub-string by removing $\s {M_{1}(S^1)_{e}}$;
 \item  [$\s{(6)}$] $\s{M_{1}(S^1)_{e}-Z(e_{0},\upsilon_1,\omega_0)-M_{2}(S^4,S^3)-C_{\eta^2}^4-M_{3}(S^1)_{e}}$;
 \item  [$\s{(7)}$] $\s{ Z(e_{0},\upsilon_0,\omega_1)-M}, ~~~M\in\{M(S^{3})_{f}, M(e_{0},e'_0), M(e_{0})_{e},M(S^{3}, S^{4}) \} $;
 \item  [$\s{(8)}$] $\s{Z(e_{0},\upsilon_0,\omega_1)\doteqdot M(S^{3},e_{0})}$;
 \item  [$\s{(9)}$] $\s{M_{1}-Z(e_{0},\upsilon_0,\omega_1)-M_{2},~~~(M_{1},M_{2})\in\{(M_{1}(e_{0})_{e},M_{2}(S^{3})_{f}), ~(M_{1}(e_{0})_{e},M_{2}(S^{3},S^{4})),~(M_{1}(e_{0},S^3),M_{2}(S^{3})_f)}$;
 \item [$\s{(10)}$]  $\s{M_{1}(S^3)_f-Z(e_{0},\upsilon_0,\omega_1)-M_{2}(e_{0},S^1)-C_{\eta^2}^4-M_{3}(S^3)_f}$ and its sub-strings 1)by removing  $\s{M_{3}(S^3)_f}$; 2) by removing  $\s{M_{1}(S^3)_f}$ and  $\s{M_{3}(S^3)_f}$;
 \item [$\s{(11)}$]  $\s{M_{1}(e_{0})_e-Z(e_{0},\upsilon_0,\omega_1)-M_{2}(S^3,S^4)-C^5}$  and its sub-string by removing $\s{M_{1}(e_{0})_e}$;
\item [$\s{(12)}$]  $\s{ Z(e_{0},\upsilon_0,\omega_0)-M}, ~~~M\in\{M(S^{1})_{e}, M(S^{4})_{f}, M(e_{0})_e, M(S^{3})_{f} \} $;
\item  [$\s{(13)}$] $\s{Z(e_{0},\upsilon_0,\omega_0)\doteqdot M, ~~M\in\{M(S^3,S^4),M(S^1,S^4),  M(S^3,e_{0}),  M(S^1,e_{0}) \}   } $;
\item [$\s{(14)}$]   $\s{M_1- Z(e_{0},\upsilon_0,\omega_0)-M_2, ~~~(M_{1},M_{2})\in\{(M_{1}(S^{3})_{f},M_{2}(S^{4})_{f}), ~(M_{1}(e_{0})_{e}, M_{2}(S^{3})_f),~(M_{1}(e_{0})_e, M_{2}(S^{4})_f),}$

      $\s{(M_{1}(S^1)_{e},M_{2}(S^{3})_{f}),~ (M_{1}(S^1)_{e},M_{2}(S^{4})_f), ~(M_{1}(e_{0})_e, M_{2}(S^{1})_{e}) \}}$;

\item [$\s{(15)}$]   $\s{M_1\doteqdot Z(e_{0},\upsilon_0,\omega_0)-M_2, ~(M_{1},M_{2})\in\{(M_{1}(e_{0},S^{1}),M_{2}(S^{3})_{f}), ~(M_{1}(e_{0},S^{1}), M_{2}(S^{3}, e_0)),~(M_{1}(e_{0},S^{1}), M_{2}(S^{4})_f),}$

  $\s{(M_{1}(e_{0},S^{1}),M_{2}(S^{4},e'_0)),~ (M_{1}(e_{0},S^{3}),M_{2}(S^{1})_e), ~(M_{1}(S^1, S^4), M_{2}(e_{0})_{e}),  (M_{1}(S^3, S^4), M_{2}(S^1)_{e}),} $

  $\s  {(M_{1}(S^3, S^4), M_{2}(S^1,S^4)),   (M_{1}(e_{0}, S^3), M_{2}(S^4)_f) \}}$;

\item [$\s{(16)}$]    $\s{M_1\doteqdot Z(e_{0},\upsilon_0,\omega_0)\doteqdot M_2, ~(M_{1},M_{2})\in\{(M_{1}(e_{0},S^{1}),M_{2}(S^{3}, S^4)), ~(M_{1}(e_{0},S^{3}), M_{2}(S^{1}, S^4))\}}$;
\item [$\s{(17)}$]  $\s{M_1\doteqdot Z(e_{0},\upsilon_0,\omega_0)<\begin{array}{c}
                                          \s M_2 \\
                                         \s M_3
                                        \end{array}
}$

for $\s{M_1=M_{1}(S^1,e_{0}), ~(M_{2}, M_{3})\in\{(M_{2}(S^3)_f, M_{3}(S^4)_f), (M_{2}(S^3)_f, M_{3}(S^4,e'_0)), (M_{2}(S^3,e_0), M_{3}(S^4)_f)\}};$

for $\s{M_1=M_{1}(S^4,S^3), ~(M_{2}, M_{3})\in\{(M_{2}(e_{0})_e, M_{3}(S^1)_e), (M_{2}(e_{0})_e, M_{3}(S^1,S^4)), (M_{2}(e_{0},S^3), M_{3}(S^1)_e)\}};$

for $\s{M_1=M_{1}(S^1,S^4), ~(M_{2}, M_{3})=(M_{2}(e_{0})_e, M_{3}(S^3)_f)};$~for $\s{M_1=M_{1}(S^3,e_{0}), ~(M_{2}, M_{3})=(M_{2}(S^1)_e, M_{3}(S^4)_f)};$

\item [$\s{(18)}$]   $\s{M_1-Z(e_{0},\upsilon_0,\omega_0)<\begin{array}{c}
                                          \s M_2 \\
                                         \s M_3
                                        \end{array}
}$ ~ where  $\s{\{M_1, M_2, M_3\}}$ is a subset of $\s {\{M(e_{0})_e, M(S^1)_e,  M(S^4)_f,   M(S^3)_f  \}};$

\item [$\s{(19)}$]  $\s{\begin{array}{c}\s M_1(e_{0})_e \\ \s M_2(S^1)_e \end{array} >Z(e_{0},\upsilon_0,\omega_0)<\begin{array}{c}
                                          \s M_3(S^3)_f \\
                                         \s M_4(S^4)_f
                                        \end{array} };$

\item [$\s{(20)}$]  $\s{M_{1}(e_{0},S^1)\doteqdot Z(e_{0},\upsilon_0,\omega_0)<\begin{array}{lc}
                                          \s {M_2(S^3)_{f}} &\\
                                         \s{ M_3(S^{4},S^1)-C_{\eta^2}^4-M_{4}(S^3)_{f}}&\\
                                        \end{array}}$
and its sub-strings 1) by removing $\s{M_{4}(S^3)_{f}}$; 2) by removing  $\s {M_2(S^3)_{f}}$ and $\s{M_{4}(S^3)_{f}};$

\item [$\s{(21)}$]   $\s{M_{1}(S^3,S^4)\doteqdot Z(e_{0},\upsilon_0,\omega_0)<\begin{array}{lc}
                                          \s {M_2(e_{0})_{e}} &\\
                                         \s{ M_3(S^1,S^{4})-C^5}&\\
                                        \end{array}}$
and its sub-string by removing $\s {M_2(e_{0})_{e}};$

\item [$\s{(22)}$]  $\s{M_{1}(S^3,S^4)\doteqdot Z(e_{0},\upsilon_0,\omega_0)<\begin{array}{lc}
                                          \s {M_2(S^1)_{e}} &\\
                                         \s{ M_3(e_{0},S^3)-C_{\eta^2}^4-M_4(S^1)_{e}}&\\
                                        \end{array}};$

\item [$\s{(23)}$]   $\s{M_{1}(S^3,S^4)\doteqdot Z(e_{0},\upsilon_0,\omega_0)<\begin{array}{lc}
                                          \s {M_2(S^1)_{e}} &\\
                                         \s{ M_3(e_{0},S^3)-C^4~}&\\
                                        \end{array}}$
 and its sub-string by removing $\s {M_2(S^1)_{e}};$
\end{itemize}

\begin{itemize}
  \item [$\s{(24)}$] $\s{Z(e_{0},S^1,\upsilon)};$
  \item [$\s{(25)}$]  $\s{Z(e_{0},S^1,\upsilon_1)-M,~M\in\{M(S^1)_e, M(S^1, e_0),  M(S^1, S^4) \}};$
  \item [$\s{(26)}$]  $\s{Z(e_{0},S^1,\upsilon_1)-M(S^1,S^4)-C^5};$
  \item [$\s{(27)}$] $\s{Z(e_{0}, S^1, \upsilon_0)-M,~M\in\{M(S^1)_e, M(S^3, S^4),  M(S^1, S^4), M(e_{0})_e, M(S^3)_{f}\}};$
   \item [$\s{(28)}$]  $\s{Z(e_{0}, S^1, \upsilon_0)\doteqdot M,~M\in\{ M(e_{0}, S^3),  M(e_{0}, S^1) \}};$
  \item [$\s{(29)}$]  $\s{M_1-Z(e_{0}, S^1, \upsilon_0)-M_2,}$
   for $\s{M_1\in\{M_{1}(S^3)_{f}, M_{1}(S^4,S^3)\},~~M_{2}\in \{M_{2}(S^1,S^4),M_{2}(S^1)_{e} \}};$ for $\s{M_{1}=M_{1}(e_{0})_{e},}$

   $\s{M_{2}\in\{M_{2}(S^1)_{e},M_{2}(S^1,S^4), M_{2}(S^3)_{f}, M_{2}(S^3,S^4)\}};$
  \item [$\s{(30)}$] $\s{M_1\doteqdot Z(e_{0}, S^1, \upsilon_0)-M_2,   (M_1,M_2)\in \{(M_{1}(S^1,e_{0}), M_{2}(S^3)_{f}),~ (M_{1}(S^1,e_{0}), M_{2}(S^3,S^4)), ~(M_{1}(S^3,e_{0}), M_{2}(S^1)) \};}$
 \item [$\s{(31)}$]  $\s{M_1-Z(e_{0}, S^1, \upsilon_0)<\begin{array}{lc}
                                          \s {M_2(e_{0})_{e}} &\\
                                         \s{ M_3}&\\
                                        \end{array},  ~M_1\in\{M_1(S^3)_f, M_{1}(S^4,S^3)\},~ M_3\in \{M_{3}(S^1)_e, M_{3}(S^1,S^4) \};}$
 \item [$\s{(32)}$]  $\s{M_{1}(S^3)_{f}-Z(e_{0}, S^1, \upsilon_0)<\begin{array}{lc}
                                          \s {M_2(e_{0})_{e}} &\\
                                         \s{ M_3(S^1,S^4)-C^5}&\\
                                        \end{array} }$  and its sub-string by removing $\s{M_2(e_{0})_{e}};$
 \item [$\s{(33)}$]   $\s{X-M_{2}(S^4,S^3)-Z(e_{0}, S^1, \upsilon_0)\doteqdot M_1(e_{0},S^1),~X=C^5~\text{or}~X=Z(e'_{0}, \omega_0)};$
 \item [$\s{(34)}$]   $\s{X-M_{2}(S^4,S^3)-Z(e_{0}, S^1, \upsilon_0)<\begin{array}{lc}
                                          \s {M_2(e_{0})_{e}} &\\
                                         \s{ M_3(S^1)_e}&\\
                                        \end{array},~X=C^5~\text{or}~X=Z(e'_{0}, \omega_0)};$
 \item [$\s{(35)}$]   $\s{M_3-Z(e'_{0}, \omega_0)-M_{2}(S^4,S^3)-Z(e_{0}, S^1, \upsilon_0)\doteqdot M_1(e_{0},S^1),~M_3\in\{ M_3(e'_{0})_e,M_3(e_0,e'_{0})\}};$
\item [$\s{(36)}$]   $\s{M_3(e'_{0})_e-Z(e'_{0}, \omega_0)-M_{2}(S^4,S^3)-Z(e_{0}, S^1, \upsilon_0)<\begin{array}{lc}
                                          \s {M_2(e_{0})_{e}} &\\
                                         \s{ M_3(S^1)_e}&\\
                                        \end{array}}$ and its sub-string by removing $\s{M_2(e_{0})_{e}};$

\item [$\s{(37)}$]   $\s{X-M_{2}(S^4,S^3)-Z(e_{0}, S^1, \upsilon_0)<\begin{array}{lc}
                                          \s {M_2(e_{0})_{e}} &\\
                                         \s{ M_3(S^1,S^4)-C^5}&\\
                                        \end{array}~ X=C^5} $ and its sub-strings 1)by removing                                   $  \s {M_2(e_{0})_{e}}$; 2)by removing $\s{X}$ and $\s{M_{2}(S^4,S^3)}$; 3)by removing  $\s{X}$, $\s{M_{2}(S^4,S^3)}$ and $  \s {M_2(e_{0})_{e}}$;

 \item [$\s{(38)}$]   $\s{X-M_{1}(S^4,S^3)-Z(e_{0}, S^1, \upsilon_0)- M_2,~M_2\in\{ M_2(e_{0})_e, M_2(S^1)_e\}, ~X=C^5~\text{or}~X= Z(e'_{0},\omega_0)};$
\end{itemize}

\begin{itemize}
  \item [$\s{(39)}$]   $\s{Z(\omega)};$
  \item [$\s{(40)}$]   $\s{Z(\omega_1)-M(S^3)};$
  \item [$\s{(41)}$]   $\s{Z(\omega_0)-M,~M\in\{M(S^1)_e, M(S^1,e_0), M(S^3)_f, M(S^3,e_0), M(S^4)_f\}};$
  \item [$\s{(42)}$]   $\s{Z(\omega_0)\doteqdot M,~M\in\{M(S^1, S^4), M(S^3,S^4)\}};$
  \item [$\s{(43)}$]   $\s{M_1-Z(\omega_0)-M_2}$, for $\s{M_1\in\{M_1(S^1)_e, M_1(e_0,S^1)\},~M_2\in \{M_2(S^3)_f, M_2(S^3, e_0), M_2(S^4)_f \}};$

   for $\s{M_1= M_1(S^4)_f }$, $\s{M_2\in \{M_2(S^3)_f, M_2(S^3, e_0) \}};$
  \item [$\s{(44)}$]   $\s{M_1-Z(\omega_0)\doteqdot M_2}$, for $\s{(M_1, M_2)\in\{(M_1(S^1), M_2(S^3, S^4)),~ (M_1(S^3)_f, M_2(S^4,S^1)),~ (M_1(e_0,S^3),  M_2(S^4,S^1)) \}};$
  \item [$\s{(45)}$]   $\s{M_1-Z(\omega_0)<\begin{array}{lc}
                                          \s {M_2} &\\
                                         \s{ M_3(S^4)_f}&\\
                                        \end{array}}$, $\s{M_1\in\{M(S^1)_e, M(e_0,S^1)\},~M_2\in \{M_2(S^3)_f, M_2(S^3, e_0) \}};$
 \item [$\s{(46)}$]  $\s{M_4(S^3)_f-Z(e_0,\upsilon_0)-M_1(e_0, S^1)-Z(\omega_0)<\begin{array}{lc}
                                          \s {M_2(S^3)_f} &\\
                                         \s{ M_3(S^4)_f}&\\
                                        \end{array}}$ and its sub-strings 1) by removing $\s{M_4(S^3)_f}$; 2) by removing $\s{ M_3(S^4)_f}$; 3) by removing $\s {M_2(S^3)_f}$ and  $\s{M_4(S^3)_f}$; 4) by removing $\s {M_2(S^3)_f}$, $\s{M_4(S^3)_f}$ and $ \s{ M_3(S^4)_f}$;
 \item [$\s{(47)}$]   $\s{M_3(S^3)_f-Z(e_0,\upsilon_0)-M_1(e_0, S^1)-Z(\omega_0)\doteqdot  M_2(S^3, S^4)}$  and its sub-string by removing $\s{M_3(S^3)_f};$
\end{itemize}
 \begin{itemize}
   \item [$\s{(48)}$]  $\s{Z(e_0,v)}$
   \item [$\s{(49)}$]  $\s{Z(e_0,v_0)-M,~M\in \{M(e_0)_e, M(e_0, e'_0), M(S^3)_f, M(S^3, S^4)\};}$
   \item [$\s{(50)}$]  $\s{Z(e_0,v_0)\doteqdot M(e_0, S^3)};$
   \item [$\s{(51)}$]  $\s{M_1-Z(e_0,v_0)-M_2,~M_1\in \{M_1(e_0)_e, M_1(e'_0, e_0)\}, M_2\in\{, M_2(S^3)_f, M_2(S^3, S^4)\};}$
 \item [$\s{(52)}$]  $\s{M_1(e_0)_e-Z(e_0,v_0)-M_2(S^3,S^4)-X,~X\in \{ Z(e'_0, \omega_0), C_{\eta}^5, C_{\eta^2}^5\}}$ and its sub-string by removing $\s{M_1(e_0)_e};$
 \item [$\s{(53)}$]  $\s{M_1(e_0)_e-Z(e_0,v_0)-M_2(S^3,S^4)-Z(e'_0, \omega_0)-M_3(e'_0)_e}$ and its sub-string by removing $\s{M_1(e_0)_e};$
 \item [$\s{(54)}$]  $\s{Z(e_0,v_0)<\begin{array}{lc}
                                          \s {M_1(e_0,e'_0)} &\\
                                         \s{ M_2(S^3,S^4)}&\\
                                        \end{array}>Z(e'_0, \omega_0)}$

 \item [$\s{(55)}$]   $\s{M_1(S^3)_f-Z(e_0,v_0)-M_2(e_0,e'_0)-Z(e'_0, \omega_0)-M_3(S^4)_f}$ and its sub-strings 1) by removing $\s{M_1(S^3)_f};$ 2) by removing $\s{M_3(S^4)_f};$

 \item [$\s{(56)}$]   $\s{M_1-Z(e_0,v_0)-M_2(e_0,S^1)-C_{\eta^2}^4,~M_1\in\{M_{1}(S^3)_f, M_{1}(S^4,S^3)\}};$
 \item [$\s{(57)}$]   $\s{M_1(e'_0,e_0)-Z(e_0,v_0)-M_2(S^3,S^4)-C^5};$
 \item [$\s{(58)}$]   $\s{C^5-M_1(S^4, S^3)-Z(e_0,v_0)-M_2(e_0,S^1)-C_{\eta^2}^4};$
  \item [$\s{(59)}$]   $\s{Z(e_0,v_0)-M_2(e_0,e'_0)-X,}$ if $\s{e'_0=S^1}$, then $\s{X=Z(e'_0, \omega_0)}$; if $\s{e'_0= C_{\eta~:n+1}^{n+3}}$, then $\s{X\in\{Z(e'_0, \omega_0), C_{\eta^2}^4\}};$
 \item [$\s{(60)}$]   $\s{M_1-Z(e_0,v_0)-M_2(e_0,S^1)-C_{\eta^2}^4-M_3,~(M_1,M_3)\in\{(M_{1}(S^3)_f, M_{3}(S^3)_f),~  (M_{1}(S^3)_f, M_{3}(S^3,S^4)),}$

      $\s{~(M_{1}(S^4,S^3), M_{3}(S^3,S^4))  \}};$

 \item [$\s{(61)}$]   $\s{M_1-Z(e_0,v_0)-M_2(e_0,S^1)-C_{\eta^2}^4-M_3(S^3,S^4)-C^5,~M_1\in\{M_{1}(S^3)_f,  M_{1}(S^4,S^3)\}};$
  \item [$\s{(62)}$]   $\s{C^5-M_1(S^4,S^3)-Z(e_0,v_0)-M_2(e_0,S^1)-C_{\eta^2}^4-M_3(S^3,S^4)-C^5};$
 \end{itemize}
\begin{itemize}
  \item [$\s{(63)}$]  $\s{Z(e'_0, \omega)};$
  \item [$\s{(64)}$]  $\s{Z(e'_0, \omega_0)-M,~M\in\{M(e'_0)_e, M(e'_0,e_0), M(S^4)_f, M(S^4,S^3) \}};$
  \item [$\s{(65)}$]  $\s{M\doteqdot M(e'_0,S^4)};$
  \item [$\s{(66)}$]  $\s{M_1-Z(e'_0, \omega_0)-M_2,~M_1\in\{M_1(e'_0)_e, M_1(e_0,e'_0)\},~ M_2\in \{M(S^4)_f, M(S^4,S^3) \}};$
  \item [$\s{(67)}$]   $\s{M_1-Z(e'_0, \omega_0)-M_2(S^4, S^3)-C^4, ~M_1\in\{M_1(e'_0)_e, M_1(e_0,e'_0)\}}$ and its sub-string by removing $\s{M_1};$
  \item [$\s{(68)}$]  $\s{M_1-Z(e'_0, \omega_0)-M_2(S^4, S^3)-C_{\eta^2}^4-M_3, ~M_1\in\{M_1(e'_0)_e, M_1(e_0,e'_0)\},~M_3\in\{M_3(S^1)_e, M_3(S^1,e_0)\}};$
\end{itemize}
\begin{itemize}
  \item [$\s{(69)}$]  $\s{C_{\eta}^4, ~~C_{\eta^2}^4, ~~C_{\eta}^5, ~~C_{\eta^2}^5};$
  \item [$\s{(70)}$]  $\s{C_{\eta^2}^4-M,~M\in\{M(S^1), M(S^3)\};~~~~~~C_{\eta}^4-M(S^3);~~~~~~C^5-M(S^4)};   $
  \item [$\s{(71)}$]  $\s{M_1(S^1)-C_{\eta^2}^4-M_2(S^3)};$
  \item [$\s{(72)}$]  $\s{M_1(S^1)-C_{\eta^2}^4-M_2(S^3,S^4)-C^5};$~~~ $\s{M_1-C_{\eta^2}^4-M_2(S^1,S^4)-C^5,~~M_1\in\{M_1(S^3)_f,~M_1(e_0, S^3)\}};$
    \item [$\s{(73)}$]  $\s{C^5-M_1(S^4,S^1)-C_{\eta^2}^4-M_2(S^3,S^4)-C^5};$
    \item [$\s{(74)}$]  $\s{C_{\eta^2}^4-M-C^5,~~M\in\{M(S^1,S^4),~M(S^3,S^4)\}}.$
\end{itemize}
For the matrix form of the polyhedra above, with row label $C_{\eta~~~:n}^{n+2}$ or $ C_{\eta^{2}~:n}^{n+3}$ ($C_{\eta~:n+1}^{n+3}$),
 $\upsilon~(\omega) \in\{1,2,3,4,6\}\subset\mathbb{Z}/12$;
$\upsilon_0~(\omega_0) \in\{3,6\}\subset\mathbb{Z}/12$;
$\upsilon_1~(\omega_1) \in\{1,2,4\}\subset\mathbb{Z}/12$.
With row label $S^n$ ($S^{n+1}$),
 $\upsilon~(\omega) \in\{1,4,6,8,9,10,12\}\subset\mathbb{Z}/24$;
$\upsilon_0~(\omega_0) \in\{6,9,12\}\subset\mathbb{Z}/24$;
$\upsilon_1~(\omega_1) \in\{1,4,8,10\}\subset\mathbb{Z}/24$;

To prove the completeness of the list, one starts from $A'_{X}$ associated to any $X$ and compute its indecomposable summands by induction on the number of indecomposable summands of $A'_{X}(2)$.

In the following we take $A'_{X}(2)$ which contains a direct summand
  $\s{(\upsilon_0\eta^2\omega_0)_{0}^0}$  $\s{=\begin{tabular}{c|c|c|}
   \multicolumn{1}{c}{}&\multicolumn{1}{c} {$\s{S^3}$}&\multicolumn{1}{c} {$\s{S^{4}}$}\\
       \cline{2-3}
   $\s{S^{0}}$ & $\s \upsilon_0$ & \s 0\\
      \cline{2-3}
     $\s{S^{1}}$ & $\s \eta^2$ & $\s \omega_0$ \\
      \cline{2-3}
     \end{tabular}}$~  as an example (this is an example of $\s{Z(e_{0},S^1,\upsilon)}$ ).

There may be entries 1 of matrix $A'_{X}(3)$  in the $S^0$-row,  $S^1$-row,  $S^3$-column,  $S^4$-column of $(\upsilon_0\eta^2\omega_0)_0^0$, denoted by $1_1, 1_2, 1_3, 1_4$ respectively; Here we consider the most complicated case: all of
  $1_1, 1_2, 1_3, 1_4$  exist and there are no entries  $\s{\begin{tabular}{c|c|}
   \multicolumn{1}{c}{}&\multicolumn{1}{c} {$\s{S^{n+4}}$}\\
       \cline{2-2}
   $\s{e'_0}$ &  \s 1\\
      \cline{2-2}
     \end{tabular}}$ and  $\s{\begin{tabular}{c|c|}
   \multicolumn{1}{c}{}&\multicolumn{1}{c} {$\s{S^{3}}$}\\
       \cline{2-2}
   $\s{e_0}$ &  \s 1\\
      \cline{2-2}
     \end{tabular}}$ of  $A'_{X}(3)$( other cases are easier, for example, if both   $\s{\begin{tabular}{c|c|}
   \multicolumn{1}{c}{}&\multicolumn{1}{c} {$\s{S^{4}}$}\\
       \cline{2-2}
   $\s{e'_0}$ &  \s 1\\
      \cline{2-2}
     \end{tabular}}$ and  $\s{\begin{tabular}{c|c|}
   \multicolumn{1}{c}{}&\multicolumn{1}{c} {$\s{S^{3}}$}\\
       \cline{2-2}
   $\s{e_0}$ &  \s 1\\
      \cline{2-2}
     \end{tabular}}$ exist, we can move the entries 1 to the place  $\s{\begin{tabular}{c|c|}
   \multicolumn{1}{c}{}&\multicolumn{1}{c} {$\s{S^{4}}$}\\
       \cline{2-2}
   $\s{S^1}$ &  $\s \omega_0$\\
      \cline{2-2}
     \end{tabular}}$ and  $\s{\begin{tabular}{c|c|}
   \multicolumn{1}{c}{}&\multicolumn{1}{c} {$\s{S^{3}}$}\\
       \cline{2-2}
   $\s{S^0}$ &  $\s {\upsilon_0}$\\
      \cline{2-2}
     \end{tabular}}$ by admissible transformations $\mathcal{G}(3)$, then $(\upsilon_1\eta^2\omega_1)_0^0$ splits out of $A'_{X}$).

$\begin{tabular}{c|c|c|cc|cc|c}
 \multicolumn{1}{c}{}&\multicolumn{1}{c}{$\s{S^{3}}$}&\multicolumn{1}{c}{$\s S^{4}$}&\multicolumn{1}{c}{$\s{\mo{r}{3}}:\s{3}$}&\multicolumn{1}{c}{$\s 4$}&\multicolumn{1}{c}{$\s{\mo{r'}{3}}:\s{3}$}&\multicolumn{1}{c}{$\s 4$}&\multicolumn{1}{c}{$\s{\cdots}$}\\
 \cline{2-8}
   $\s S^0$& $\s\upsilon_0$& \s 0 & $\s 1_1$&\s 0 &\s 0& ~~\s 0&\s 0\\
      \cline{2-8}
   $\s S^{1}$& $\s\eta^2$& $\s \omega_0$ &\s  0& \s 0 &\s 0& ~~ $\s 1_2$&\s 0\\
      \cline{2-8}
   $\s{\mo{s}{0}}:\s{0}$&$\s 1_3$& \s 0&&&&  ~~&\\
    $\s{1}$& \s 0 &$\s 0$ & & &&~ &\\
     \cline{2-8}
      $\s{\mo{s'}{0}}:\s{0}$&\s 0 & \s 0 & &&&  &\\
    $\s{1}$& \s 0 &$\s 1_4$ & & &&~& \\
    \cline{2-8}
   &\s 0&\s 0 && &&& \\
 \end{tabular}$~

  Assume $A'_{X}(3)=\bigoplus A'_{j}(3)$, $A'_{j}(3)$ is a string object of  $\mathrm{ind} \mathscr{A}(3)$ for any $j$.\newline
  \newline
 \textbf{Claim 1} $M(e_0, S^3)$ and $M(e'_0, S^4)$   (resp. $M(e_0, e'_0)$ and $M( S^3, S^4)$ ) are direct summands of  $A'_{X}(3)$ simultaneously, then  $\s{M(e_0, S^3)\doteqdot (\upsilon_0\eta^2\omega_0)_0^0 \doteqdot M(e'_0, S^4)}$ (resp. $\s{M(S^1, S^0)\doteqdot (\upsilon_0\eta^2\omega_0)_0^0\doteqdot M(S^3, S^4)}$) splits out.
\begin{proof}[proof of Claim 1]
  If $A'_{j_1}(3)=M(e_0, S^3)$ and  $A'_{j_2}(3)=M(e'_0, S^4)$, then  $\s{M(e_0, S^3)\doteqdot (\upsilon_0\eta^2\omega_0)_0^0}$ $\s{\doteqdot M(e'_0, S^4)}$ splits out since we can move 1 in $e_0$-row of $A'_{j_1}(3)$ to  $S^0$-row  of $(\upsilon_0\eta^2\omega_0)_0^0 $ and move 1 in $e'_0$-row of $A'_{j_2}(3)$ to  $S^1$-row of $(\upsilon_0\eta^2\omega_0)_0^0 $ by admissible transformations $\mathcal{G}(3)$; Similarly, if $A'_{j_1}(3)=M(e_0, e'_0)$ and   $A'_{j_2}(3)=M( S^3, S^4)$, then $\s{M(S^1, S^0)\doteqdot (\upsilon_0\eta^2\omega_0)_0^0 \doteqdot M(S^3, S^4)}$ splits out of $A'_X$;
\end{proof}

 From \textbf{Claim 1}, we only need to consider the following cases
 \begin{itemize}
   \item [(i)] $\s{ M_1(S^0, S^1)\doteqdot (\upsilon_0\eta^2\omega_0)_0^0 <\begin{array}{lc}
                                          \s {M_2(S^3)} &\\
                                         \s{ M_3(S^{4})}&\\
                                        \end{array}}$, $ \s {M_2(S^3)\in \{M_2(S^3)_f, M_2(S^3, e_0)\}}$, $ \s {M_3(S^4)\in \{M_3(S^4)_f, M_3(S^4, e'_0)\}}$

 ($\s{ M_2(S^3, e_0)}$ and $\s{ M_3(S^4, e'_0)}$ do not appear simultaneously);
   \item [(ii)] $\s{ M_1(S^3, S^4)\doteqdot (\upsilon_0\eta^2\omega_0)_0^0 <\begin{array}{lc}
                                          \s {M_2(S^0)} &\\
                                         \s{ M_3(S^{1})}&\\
                                        \end{array}}$, $ \s {M_2(S^0)\in \{M_2(S^0)_e, M_2(S^0, S^3)\}}$, $ \s {M_3(S^1)\in \{M_3(S^1)_e, M_3(S^1, S^4)\}}$;
   \item [(iii)]$\s{ M_1(S^1, S^4)\doteqdot (\upsilon_0\eta^2\omega_0)_0^0 <\begin{array}{lc}
                                          \s {M_2(S^0)} &\\
                                         \s{ M_3(S^{3})}&\\
                                        \end{array}}$, $ \s {M_2(S^0)=M_2(S^0)_e}$, $ \s {M_3(S^3)\in \{M_3(S^3)_f, M_3(S^3, e_0)\}}$

   (note that if $\s{M_2(S^0)=M_2(S^0, e'_0)}$, then it will be the case (i) by admissible transformations);
   \item[(iv)]$\s{ M_1(S^0, S^3)\doteqdot (\upsilon_0\eta^2\omega_0)_0^0 <\begin{array}{lc}
                                          \s {M_2(S^1)_e} &\\
                                         \s{ M_3(S^{4})_f}&\\
                                        \end{array}}$(this can split out of  $A'_X$, i.e. in (17) of \textbf{List*});
   \item[(v)]$\s{ \begin{array}{lc}
                                          \s {M_2(S^0)_e} &\\
                                         \s{ M_3(S^{3})_f}&\\
                                        \end{array}> (\upsilon_0\eta^2\omega_0)_0^0 <\begin{array}{lc}
                                          \s {M_2(S^1)_e} &\\
                                         \s{ M_3(S^{4})_f}&\\
                                        \end{array}}$(this can split out of  $A'_X$,  i.e. in (19) of \textbf{List*}).
 \end{itemize}
By induction it suffices to  find the indecomposable strings  containing $(\upsilon_0\eta^2\omega_0)_0^0$ for the case (i) above, and the other cases are similar.

If a matrix $X\in \mathrm{ind}\mathscr{A}(2)$ contains $x_1$-row, $\cdots$, $x_k$-row and $y_1$-column, $\cdots$, $y_t$-column, then we denote it by $X(x_1,\cdots, x_k, y_1,\cdots, y_t)$, where $x_i\in e_0\coprod e'_0$, $y\in\{S^3,S^4\}$ (hence $k,t\leq 2$).

 (i1)~~If $\s {M_2(S^3)=M_2(S^3)_f}$, $\s {M_3(S^4)=M_3(S^4)_f}$, then $\s{ M_1(S^0, S^1)\doteqdot (\upsilon_0\eta^2\omega_0)_0^0 <\begin{array}{lc}
                                          \s {M_2(S^3)_f} &\\
                                         \s{ M_3(S^{4})_f}&\\
                                        \end{array}}$ splits out of $A'_X$.

 (i2)~~If $\s {M_2(S^3)=M_2(S^3)_f}$, $\s {M_3(S^4)=M_3(S^4, e'_0)}$, then  $\s{ M_1(S^0, S^1)\doteqdot (\upsilon_0\eta^2\omega_0)_0^0 <\begin{array}{lc}
                                          \s {M_2(S^3)_f} &\\
                                         \s{ M_3(S^{4},e'_0)}&\\
                                        \end{array}}$ splits out of $A'_X$ or else we have the  string $(\clubsuit):$ $\s{ M_1(S^0, S^1)\doteqdot (\upsilon_0\eta^2\omega_0)_0^0 <\begin{array}{lc}
                                          \s {M_2(S^3)_f} &\\
                                         \s{ M_3(S^{4},e'_0)-X(e'_0)}&\\
                                        \end{array}}$.
 \newline
 \textbf {Claim 2} If string $(\clubsuit)$ is a substring of  indecomposable string $X$, then $X(e'_0)$ dose not contain $S^3$-column, i.e.  $X(e'_0)\neq X(e_0, e'_0,S^3,S^4)$.
\begin{proof}[proof of Claim 2]
Suppose that $X(e'_0)= X(e_0, e'_0,S^3,S^4)$ in the later case of (i2).
Note that $(\upsilon_0\eta^2\omega_0)_0^0$ is connected to $M_3(S^4, e'_0)$ by $S^4$, thus $M_3(S^4, e'_0)$ is connected to $X(e_0, e'_0,$ $S^3,S^4)$ by $e'_0$.
 $X(e_0, e'_0,S^3,S^4)$ maybe also connected to $M_4(S^3)$ by $S^3$, to $M_5(e_0)$ by $e_0$ and to $M_6(S^4)$ by $S^4$.

  By the \textbf{Claim 1}, $M_4(S^3)=M_4(S^3)_f$. $M_5(e_0)=M_5(e_0)_e$ or $M_5(e_0,e'_0)$. $M_6(S^4)=M_6(S^4)_f$ or $M_6(S^4,e'_0)$. For $M_5(e_0)=M_5(e_0,e'_0)$ and $M_6(S^4)=M_6(S^4,e'_0)$ we will show that the following ``closed" string $(\maltese)$ can not be a sub-string of this indecomposable string  $X$, then by exchanging  the nonzero entry $1_4$ in the $S^4$-column of $(\upsilon_0\eta^2\omega_0)_0^0$  and nonzero entry in the $S^4$-column of $X(e_0, e'_0,S^3,S^4)$ by admissible transformations $\mathcal{G}(3)$ we find that $X(e_0, e'_0,S^3,S^4)$ and $(\upsilon_0\eta^2\omega_0)_0^0 $ are not in the same indecomposable string.
$$
 \footnotesize{\xymatrix@R=0.1cm{
                &        \s{ M_6(S^4,e'_0)} \ar@{-}[r]  &\s{X'_1} \ar@{-}[r] &  \s{M'_1} \ar@{-}[r] &\s{X'_2}\\
  \s{ X(e_0, e'_0,S^3,S^4)}~~~~~~~\ar@{-}[ur] \ar@{-}[dr]        &     &  && &\ast\ar@{.}[ul] \ar@{.}[dl]& ~~ ~~~~~(\maltese) \\
                &        \s{ M_5(e_0,e'_0)}   \ar@{-}[r]     &    X'_n  \ar@{-}[r] & \s{M'_{n-1}} \ar@{-}[r]&\s{X'_{n-1}} }}
 $$
                            $$ \text{where}~ M'_{i}~\text{  are string objects of~} \mathrm{ind}\mathscr{A}(3), X'_{i}\in \mathrm{ind}\mathscr{A}(2).$$

Prove it by contradiction. Suppose that the above ``closed" string $(\maltese)$ is a sub-string of  indecomposable string  $X$.  By \textbf{Claim 1} $M'_i\neq M'_i(S^4,S^3)$ and $M'_i\neq M'_i(e_0, S^3)$.  So $S^3$-column can not appear in $M'_i$. From $M_6(S^4, e'_0)$ we know that $X'_1$ contains a $e'_0$-row, so the sub-string  $\s{X'_1-M'_1}$ is $\s{X'_1(e'_0,S^4)-M'_1(S^4)}$ or  $\s{X'_1(e_0, e'_0,S^3,S^4)-M'_1(e_0)}$, which implies that $M'_1=M'_1(S^4,e'_0)$  or $M'_1(e_0,e'_0)$. Hence $X'_2$ contains a $e'_0$-row. By the same analysis, we get  $M'_2=M'_2(S^4,e'_0)$  or $M'_2(e_0,e'_0)$. Keep going, we will get
$M'_i=M'_i(S^4,e'_0)$  or $M'_i(e_0,e'_0)$, and $M'_i$ is connected to $X'_i$ by $e'_0$ ($i=1,2,\cdots, n-1$). Specially,
$M'_{n-1}$ is connected to $X'_n$  by $e'_0$. However $M_5(e_0,e'_0)$ is also connected to $X'_n$  by $e'_0$, we get a contradiction.

\end{proof}

From \textbf{Claim 2}, there are following two possibilities for $X(e'_0)$ in the string $(\clubsuit)$.

  $\bullet$~~If $X(e'_0)$ contains $S^4$-column but not $S^3$-column, i.e. $X(e'_0)=X(e'_0, S^4)$,  then it is easy to observe that $X(e'_0, S^4)$ and $(\upsilon_0\eta^2\omega_0)_0^0 $ are not in the same indecomposable string by moving $1_4$ to the $S^4$-column of $X(e'_0, S^4)$.

 $\bullet$~~If $X(e'_0)$ contains $S^3$-column but not $S^4$-column, i.e. $X(e'_0)=X(e'_0, S^3)$, then $e'_0$ in $M_3(S^{4},e'_0)$ must be $S^1$ and $X(e'_0, S^3)=C_{\eta^2}^4$. So
  $\s{ M_1(S^0, S^1)\doteqdot(\upsilon_0\eta^2\omega_0)_0^0 <\begin{array}{lc}
                                          \s {M_2(S^3)_f} &\\
                                         \s{ M_3(S^{4},S^1)-C_{\eta^2}^4}&\\
                                        \end{array}}$ splits out of $A'_X$ or else $\s{ M_1(S^0, S^1)\doteqdot(\upsilon_0\eta^2\omega_0)_0^0 <\begin{array}{lc}
                                          \s {M_2(S^3)_f} &\\
                                         \s{ M_3(S^{4},S^1)-C_{\eta^2}^4}-M_4(S^3)&\\
                                        \end{array}}$.
  In the later case, $M_4(S^3)=M_4(S^3)_f$, and this string is indecomposable and it is in (20) of the \textbf{List*}.¡¡

(i3)~~~If $\s {M_2(S^3)=M_2(S^3,e_0)}$, $\s {M_3(S^4)=M_3(S^4)_f}$, then $\s{ M_1(S^0, S^1)\doteqdot (\upsilon_0\eta^2\omega_0)_0^0 <\begin{array}{lc}
                                          \s {M_2(S^3,e_0)} &\\
                                         \s{ M_3(S^{4})_f}&\\
                                        \end{array}}$ splits out of $A'_X$ or else $\s{ M_1(S^0, S^1)\doteqdot (\upsilon_0\eta^2\omega_0)_0^0 <\begin{array}{lc}
                                          \s {M_2(S^3,e_0)-X(e_0, S^3)} &\\
                                         \s{ M_3(S^{4})_f}&\\
                                        \end{array}}$. But in the later case,  $X(e_0, S^3)$  and  $(\upsilon_0\eta^2\omega_0)_0^0$ are not in the same indecomposable string. The proof is similar to that of case (i2).

 In conclusion, we get all the indecomposable congruence classes of $\mathcal{A}\dag_{n+3} \mathcal{B}$:

 1) $S^{n+k}, 0\leq k\leq 5$, $C_{\eta}^{n+2}$, $C_{\eta}^{n+3}$, $C_{\eta^2}^{n+3}$; 2) $\mo{r}{n+k}, k=0,1,4, r\in \mathbb{N}_+;$ 3) band objects $B(w,z,\pi)$ and string objects in \textbf{List} $\mathrm{ind} \mathscr{A}(3)$; 4) string objects in \textbf{List*}.

For $X\in \mathbf{F}_{n(2)}^5$, if $H_{n+3}X$ is not 3-torsion free, then the matrix $A_{X}$ realizing $X$ is also a block matrix with block $\gamma_{ij}$ which has entries from the (ij)-th cell of Table $\Gamma(A,B)$ except that the  $\s{\begin{tabular}{c|c|}
   \multicolumn{1}{c}{}&\multicolumn{1}{c} {$\s{S^{n+3}}$}\\
       \cline{2-2}
   $\s{S^{n+3}}$ &$\mathbb{Z}$\\
      \cline{2-2}
     \end{tabular}}$~~--block of $A_X$ is nonzero.  Since $X$ is 2-torsion free and $M_{3^{r}q}^{n+3}\simeq M_{3^{r}}^{n+3}\vee M_{q}^{n+3}$ ( $M_{q}^{n+3}$ can split out of the cofiber where odd number $q$ is not divided by 3 ), $\s{\begin{tabular}{c|c|}
   \multicolumn{1}{c}{}&\multicolumn{1}{c} {$\s{S^{n+3}}$}\\
       \cline{2-2}
   $\s{S^{n+3}}$ &$\mathbb{Z}$\\
      \cline{2-2}
     \end{tabular}}$~~--block of $A_X$  can be diagonalized to $diag(\lambda_1, \lambda_2, \cdots, \lambda_t, 0,\cdots, 0 )$ with $\lambda_{i}=3^{r_{i}}$ where $r_{i}$ is a nonnegative integer for $i=1,2,\cdots, t$.

 We can use nonzero entry $\lambda_{i}$ to make nonzero entries in $\mathbb{Z}/2$-blocks of the same row and the same column of $A_X$ zero. Moreover, we can also use $\lambda_{i}$ to eliminate the $2$-primary component of  entries in $\mathbb{Z}/24$-blocks or $\mathbb{Z}/12$-blocks of the same column of $A_X$.
 Hence, $A_X$ contains the following sub-matrix under admissible transformations
 $$\s{\begin{tabular}{c|c|}
   \multicolumn{1}{c}{}&\multicolumn{1}{c} {$\s{S^{n+3}}$}\\
       \cline{2-2}
   $\s{S^{n}}$ &\s 8\\
      \cline{2-2}
       $\s{S^{n+3}}$ &$\s{3^{r}}$\\
      \cline{2-2}
     \multicolumn{1}{c}{}& \multicolumn{1}{c}{$\s{C_{8}^{n+4}(r)}$}
     \end{tabular}}~;~~\s{\begin{tabular}{c|c|}
   \multicolumn{1}{c}{}&\multicolumn{1}{c} {$\s{S^{n+3}}$}\\
    \cline{2-2}
       $\s{S^{n+3}}$ &$\s{3^{r}}$\\
       \cline{2-2}
    $\s{C_{\eta}^{n+2}}:\s{n}$ & $\s 4$\\
      \quad ${\s{n+2}}$ & \s 0 \\
      \cline{2-2}
       \multicolumn{1}{c}{}& \multicolumn{1}{c}{$\s{(\eta4)_{0}^{1}(r)}$}
     \end{tabular}}~;~~\s{\begin{tabular}{c|c|}
   \multicolumn{1}{c}{}&\multicolumn{1}{c} {$\s{S^{n+3}}$}\\
    \cline{2-2}
       $\s{S^{n+3}}$ &$\s{3^{r}}$\\
       \cline{2-2}
    $\s{C_{\eta^2}^{n+3}}:\s{n}$ & $\s 4$\\
      \quad ${\s{n+3}}$ & \s 0 \\
      \cline{2-2}
       \multicolumn{1}{c}{}& \multicolumn{1}{c}{$\s{(\eta^24)_{0}^{1}(r)}$}
     \end{tabular}}~;~~\s{\begin{tabular}{c|c|}
   \multicolumn{1}{c}{}&\multicolumn{1}{c} {$\s{S^{n+3}}$}\\
    \cline{2-2}
       $\s{S^{n+3}}$ &$\s{3^{r}}$\\
       \cline{2-2}
    $\s{M_{3^s}^{n}}:\s{n}$ & $\s 1$\\
      \quad ${\s{n+1}}$ & \s 0 \\
      \cline{2-2}
       \multicolumn{1}{c}{}& \multicolumn{1}{c}{$\s{M(-,S^3)(r)}$}
     \end{tabular}}$$
     where $r\in \mathbb{N}_+$.

  Note that if $A_X$ contains $C_{8}^{n+4}(r)$ (resp. $(\eta4)_{0}^{1}(r)$, $(\eta^24)_{0}^{1}(r)$), then $C_{8}^{n+4}(r)$ (resp. $(\eta4)_{0}^{1}(r)$, $(\eta^24)_{0}^{1}(r)$) must be a direct summand of $A_X$. While all indecomposable polyhedra containing $M(-,S^3)(r)$ can be obtained by replacing $M(-,S^3)$  in  string objects $X$ of \textbf{List*} with $M(-,S^3)(r)$ where $M(-,S^3)$ does not connect to any element of  $\mathrm{ind}\mathscr{A}(2)$ in $X$ by $S^{n+3}$.

  Therefore, we can get all the indecomposable polyhedra $X\in \mathrm{ind}\mathbf{F}_{n(2)}$ with $H_{n+3}X$  not $3$-torsion free, which are denoted by \textbf{List**}.

   \bigskip
 \textbf{List**:}

 \begin{itemize}
   \item $M_{3^r}^{n+3}$, $C_{8}^{n+4}(r)$,  $(\eta4)_{0}^{1}(r)$, $(\eta^24)_{0}^{1}(r)$;
   \item all indecomposable polyhedra obtained by replacing $M(-,S^3)$  in  string objects $X$ of \textbf{List*} with $M(-,S^3)(r)$ where $M(-,S^3)$ does not connect to any element of  $\mathrm{ind}\mathscr{A}(2)$ in $X$ by $S^{n+3}$;
 \end{itemize}

   \bigskip

Combining with Lemma \ref{lemmaindFn25} and Lemma  \ref{lemmaindA+B2tor}, we complete the classification of indecomposable congruence classes of $\mathbf{F}_{n(2)}^5$.

\begin{theorem}[\textbf{Main theorem}]\label{maintheorem}
The indecomposable congruence classes of $\mathbf{F}_{n(2)}^5$ are as follows:
\begin{itemize}
 \item [(i)] $S^{n+k}, 0\leq k\leq 5$,  $C_{\eta}^{n+2}$, $C_{\eta}^{n+3}$, $C_{\eta^2}^{n+3}$;
  \item [(ii)] $M_{p^r}^{n+k},  0\leq k\leq 4$, prime $p\neq 2$, $r\in \mathbb{N}_+$;
   \item [(iii)] band objects $B(w,z,\pi)$ and string objects in \textbf{List} $\mathrm{ind}\mathscr{A}(3)$;
   \item[(iv)] string objects in \textbf{List*} and \textbf{List**}.

   \end{itemize}
\end{theorem}


\begin{thebibliography}{9}

\bibitem{RefTF5cels}
Baues H J, Drozd Y A. The homotopy classification of (n-1)-connected (n+4)-dimensional polyhedra with torsion free homology, $n\geq5$. Expositiones Mathematicae, 1999, 17: 161-180
\bibitem{RefTF6cels}
Baues H J, Drozd Y A.  Classification of stable homotopy types  with torsion-free homology. Topology, 2011, 40: 789-821
\bibitem{RefBD}
Baues H J, Drozd Y A.  Indecomposable homotopy types with at most two non-trivial homology groups. In: Groups of Homotopy Self-Equivalences and Related Topics. Contemporary Mathematics, 2001, 274: 39¨C56
\bibitem{RefBH}
Baues H J, Hennes M. The homotopy classification of (n-1)-connected (n+3)-dimensional polyhedra, $n\geq 4$. Topology, 1991, 30: 373-408
\bibitem{RefBunchChains}
Bondarenko V M, Representations of bundles of semichained sets and their applications. St Petersburg Mathematical Journal, 1991, 3 (5): 38-61
\bibitem{RefChang}
Chang S C. Homology invariants and continuous mappings. Proc. Roy. Soc. London. ser. A, 1950, 202: 253-263
\bibitem{RefCo}
Cohen J M. Stable homotopy. Lect. Notes Math. vol. 165. Berlin Heidelberg New York: Springer-Verlag, 1970
\bibitem{RefFinitquadra}
Drozd Y A. Finitely generated quadratic modules. Manuscripta Mathematica. 2001, 104(2): 239-256
\bibitem{RefDrMSP}
Drozd Y A. Matrix problems and stable homotopy types of polyhedra. Central European J. Math, 2004, 2: 420-447
\bibitem{RefDrFP}
Drozd Y A. On classification of torsion free polyhedra. Preprint series, Max-Planck-Institut f\"{u}r Mathematik (Bonn), 2005, 92
\bibitem{RefDrMTS}
Drozd Y A. Matrix problems, triangulated categories and stable homotopy types. Sao Paulo Journal of Mathematical Sciences, 2010, 4: 209-249
\bibitem{PZ}
Pan J Z, Zhu Z J. The classification of 2 and 3 torsion free polyhedra.  Acta Mathematica Sinica, English Series, 2015, 31.11: 1659-1682
\bibitem{PZ2}
Pan J Z, Zhu Z J. Stable homotopy classification of $A_{n}^4$-polyhedra with 2-torsion free homology. Science China Mathematics, 2016, 59.6: 1141-1162
\bibitem{RefRMS}
Switzer R M. Algebraic Topology-Homology and Homotopy. Berlin: Springer-Verlag, 1975
\bibitem{RefUnsold}
Uns\"{o}ld H M. $A_{n}^4$-polyhedra with free homology. Manuscripta mathematica, 1989, 65: 123-146





\end{thebibliography}
\end{document}